\documentclass[a4paper,12pt]{article}
\usepackage[french, ngerman, italian, english]{babel}
\usepackage{amsmath,amsthm, amssymb}
\usepackage{setspace}
\usepackage{anysize}
\usepackage{url}
\usepackage{graphicx}
\usepackage[colorlinks=true]{hyperref}
\usepackage{mathptmx}
\usepackage{paralist}
\usepackage{verbatim}
\usepackage[utf8]{inputenc}

\theoremstyle{plain}
\newtheorem{thm}{\bf Theorem}[section]
\newtheorem{prop}[thm]{\bf Proposition}
\newtheorem{lem}[thm]{\bf Lemma}
\newtheorem{cor}[thm]{\bf Corollary}

\theoremstyle{definition}
\newtheorem{defn}{\bf Definition}[section]
\theoremstyle{remark}

\newtheorem{es}[thm]{\bf Example}
\newtheorem{ess}[thm]{\bf Examples}
\newtheorem{os}[thm]{Remark}

\def \tt{\tau}
\def \aa{{\bf a}}
\def \bb{{\bf b}}
\def \cc{{\bf c}}
\def \dd{{\bf d}}
\def \id{\operatorname{id}}
\def \e{\operatorname{e}}
\def \NN{\mathbb N}
\def \ZZ{\mathbb Z}

\def \KK{\mathbb K}
\def \mm{\mathfrak{m}}
\def \vv{\mathfrak{v}}
\def \jj{\mathfrak{j}}

\def \AA{A(H,\omega)}

\def \depth{\operatorname{depth}}
\def \Proj{\operatorname{Proj}}
\def \gdim{\operatorname{gdim}}
\def \ara{\operatorname{ara}}
\def \cd{\operatorname{cd}}
\def \projdim{\operatorname{projdim}}
\def \rank{\operatorname{rank}}
\def \AG{\bar{A}(G)}
\def\a{\ensuremath{\alpha}}

\begin{document}
\title{Dimension, depth and zero-divisors \\of the algebra of basic $k$-covers of a graph.}
\author{Bruno Benedetti\thanks{%
Supported by DFG via the Berlin Mathematical School; partially supported by PRAGMATIC 08.} \\
\footnotesize Inst.\ Mathematics, MA 6-2\\
\footnotesize TU Berlin\\
\footnotesize D-10623 Berlin, Germany\\
\footnotesize \url{benedetti@math.tu-berlin.de}
\and \setcounter{footnote}{8} \!\!\!\!Alex Constantinescu%
\thanks{%
Supported by Phd of Universit\'a di Genova; partially supported by PRAGMATIC 08.} \\
\footnotesize Dipartimento di Matematica\\
\footnotesize Universit\'{a}  di Genova\\
\footnotesize I-16146 Genoa, Italy\\
\footnotesize \url{constant@dima.unige.it}
\and \setcounter{footnote}{3} Matteo Varbaro%
\thanks{%
Supported by Phd of Universit\'a di Genova; partially supported by PRAGMATIC 08.} \\
\footnotesize Dipartimento di Matematica\\
\footnotesize Universit\'{a}  di Genova\\
\footnotesize I-16146 Genoa, Italy\\
\footnotesize \url{varbaro@unige.it}}
\date{{\small \today}}
\maketitle

\begin{abstract}
\noindent We study the basic $k$-covers of a bipartite graph $G$; the algebra $\AG$ they span, first studied by Herzog, is the fiber cone of the Alexander dual of the edge ideal. We characterize when $\AG$ is a domain in terms of the combinatorics of $G$; if follows from a result of Hochster that when $\AG$ is a domain, it is also Cohen-Macaulay. We then study the dimension of $\AG$ by introducing a geometric invariant of bipartite graphs, the ``graphical dimension''. We show that the graphical dimension of $G$ is not larger than $\dim(\AG)$, and equality holds in many cases (e.g. when $G$ is a tree, or a cycle). Finally, we discuss applications of this theory to the arithmetical rank.
\end{abstract}

\section*{Introduction}

To each homogenous ideal $I$ of a standard graded $\KK$ algebra $S$, we can associate a second standard graded $\KK$-algebra, the so-called \textit{fiber cone} of $I$, defined as
\[F_{\mm}(I):=\bigoplus_{k\geq 0}I^k/\mm I^k,\]
where $\mm= \oplus_{k>0}S_k$ is the maximal irrelevant ideal of $S$.

Fiber cones have been introduced by Northcott and Rees \cite{NR} and extensively studied in literature ever since. Some typical questions (in general hard to solve) are the following: given an ideal $I$,
\begin{compactitem}[$\qquad \rhd$]
 \item is $F_{\mm}(I)$ a domain? is it Cohen-Macaulay?
\item what is the Krull dimension of $F_{\mm}(I)$? (This invariant is known as the ``analytic spread of $I$''.)
\end{compactitem}

At the summer school ``Pragmatic 2008'', with Herzog and Welker as teachers, we studied the special case where $S$ is the polynomial ring in $n$ variables, and $I$ is a monomial ideal associated to a bipartite graph $G$ (more precisely, $I$ is the Alexander dual of the edge ideal of $G$). We will denote by $\bar{A}(G)$ the fiber cone of this ideal $I$. Consistently with the general picture, we focused on the following questions:

\begin{compactitem}[$\qquad \rhd$]
\item When is $\bar{A}(G)$ a domain? When is it Cohen-Macaulay?
\item What is the dimension of $\bar{A}(G)$?
\end{compactitem}

Here is a summary of some of the results contained in the present paper.

\vspace{2mm}

We give a simple, purely combinatorial characterization of those bipartite graphs $G$ for which $\bar{A}(G)$ is a domain ({\bf Theorem \ref{thm:domain}}). Namely, for \textit{every} non-isolated vertex $i$ of $G$ there must be an edge $\{i,j\}$ in $G$ such that
\[ \hbox{ If } \{i,i'\} \hbox{ is an edge of $G$ and } \{j,j'\} \hbox{ is an edge of $G$, then } \{i',j'\} \hbox{ is also an edge of $G$.}\]
In case $G$ satisfies this condition, we show via a theorem of Hochster that $\bar{A}(G)$ is actually a Cohen-Macaulay normal domain ({\bf Proposition \ref{thm:graziehochster}}). Moreover we prove in {\bf Theorem \ref{hibiring}} that $\AG$ is an Hibi ring provided that $G$ is unmixed (in this case one can easily see that $\AG$ is a domain). This has some interesting consequences: for example it allows us to characterize all Cohen-Macaulay bipartite graphs for which $\AG$ is Gorenstein ({\bf Corollary \ref{gor}}).

When $\bar{A}(G)$ is not a domain, it may be non-Cohen-Macaulay as well. However, we show that $\depth(\bar{A}(G))\geq 2$ for all bipartite graphs $G$: see {\bf Theorem \ref{thm:depth}}. This fact has some interesting consequences, among which include that the projective scheme $\Proj(\bar{A}(G))$ is always connected. Anyway we show an example of a bipartite graph $G$ for which $\AG$ is not equidimensional, and thus not Cohen-Macaulay ({\bf Example \ref{G not CM}}).

\vspace{2mm}

Furthermore: we prove that the dimension of $\bar{A}(G)$ is $ \frac{n}{2}$ when $G$ is a cycle ({\bf Proposition \ref{thm:dimcycle}}), and $\lfloor\frac{n}{2} \rfloor + 1$ when $G$ is a path ({\bf Corollary \ref{thm:dimpath}}) or a Cohen Macaulay graph ({\bf Corollary \ref{thm:dimCM}}). These results are not trivial, as the ``obvious'' approach to the problem -- that is, to asymptotically count basic $k$-covers\footnote{See the next paragraph for the definition.} -- leads to long, nasty calculations. We followed in fact another road, introducing a numerical invariant -- the ``graphical dimension $\gdim(G)$'' -- for each bipartite graph $G$.

For a generic bipartite graph $G$, we show that
\[\gdim (G) \leq \dim \bar{A}(G) \ \ \ \mbox{({\bf Theorem \ref{thm:inequalities}})}.\]
Indeed, for trees and cycles, one has equality:
\[\gdim (G) = \dim \bar{A}(G) \ \ \ \mbox{({\bf Theorem \ref{6}} and {\bf Proposition \ref{thm:dimcycle}})}.\]

We do not know whether equality holds true in general, or not.

Finally, in the last section, we discuss:
\begin{compactitem}
\item the dimension issue in a more general setting (that considers hypergraphs instead of bipartite graphs);
\item unmixed graphs and Hibi rings;
\item applications to the arithmetical rank of the Alexander dual of the edge ideal of $G$.
\end{compactitem}



\subsection{Definitions and notation}
\subsubsection{Bipartite graphs and vertex covers}
Let $G$ be a simple graph on $n$ vertices, labeled by $1$, $\ldots, n$. A vector $\aa=(a_1, \ldots, a_n)$ a \textit{$k$-cover} of $G$ if:
\begin{compactitem}
 \item $k$ and all $a_i$'s are nonnegative integers;
\item not all $a_i$'s are zero;
\item for every edge $\{i,j\}$ of $G$ one has $a_i + a_j \geq k$.
\end{compactitem}

This concept generalizes the graph theoretical notion of \textit{vertex cover}, which is a subset of $[n]$ that has non-empty intersection with all the edges of $G$. Vertex covers in our language are just 1-covers whose entries are either $0$ or $1$.

Two covers can be ``summed'' vertex-wise: a $k$-cover plus a $k'$-cover gives a $(k+k')$-cover. A $k$-cover $\aa$ is said to be \textit{decomposable} if $\aa=\bb+\cc$, for a suitable $h$-cover $\bb$ and a suitable $(k-h)$-cover $\cc$. A $k$-cover $\aa$ is said to be  \emph{non-basic} if it can be decomposed into a $k$-cover $\aa'$, and a $0$-cover $\aa''$. \textit{Indecomposable} (resp.  \textit{basic}) is of course contrary to decomposable (resp. non-basic). We also say that a $k$-cover $\aa$ can be ``lopped at $i$'' if replacing $a_i$ with $a_{i}-1$ in the vector $\aa$, gives again a $k$-cover.

\begin{es}
Let $G$ be a hexagon, with the vertices labelled clockwise. Then, $\aa=(1,0,1,1,0,1)$ is a basic $1$-cover;  $\aa=(1,1,0,1,1,0)$ is also a basic $1$-cover; yet their sum $\cc=\aa+\bb=(2,1,1,2,1,1)$ is a $2$-cover that is not basic ($c$ can be lopped at $1$). This shows that $\bar{A}(G)$ is not a domain when $G$ is a hexagon: $\aa$ and $\bb$ are both zero-divisors.

The algebra has however depth at least $2$, since we will see that  $\dd=(1,0,1,0,1,0)$ and ${\bf e}=(0,1,0,1,0,1)$ are non-zero-divisors and form a regular sequence. As a consequence of Theorem \ref{thm:dimcycle} we will see that the dimension of $\bar{A}(G)$  is 3.
\end{es}

Recall that a graph $G$ is \textit{bipartite} if its vertex set $[n]:=\{1, \ldots, n\}$ can be partitioned into subsets $A$, $B$ such that all the edges of $G$ link one vertex in $A$ with one vertex in $B$. In case every vertex in $A$ is adjacent to every vertex of $B$, we say that $G$ is \textit{bipartite complete}, and we write $G = K_{a, b}$, where $a = |A|$ and $b=|B|$. A graph $G$ is bipartite if and only if  $G$ does not contain an odd cycle. Moreover, the disjoint union of $k$ graphs $G_1, \ldots, G_k$ is bipartite if and only if each $G_i$ is bipartite. 
Finally, it was established in \cite[Theorem 5.1]{HHT}) that, providing the graph $G$ has at least one edge:
\begin{compactenum}
\item If $\aa$ is a $k$-cover with $k \geq 3$, $\aa$ is decomposable;
\item the graph $G$ is bipartite if and only if $G$ has no indecomposable $2$-covers.
\end{compactenum}
Notice that a bipartite graph might admit a basic $k$-cover even when $k \geq 2$: for example, if $G$ consists in two vertices and a single edge, and $\aa$ is the vector $(0,k)$, clearly $\aa$ is decomposable as $\aa= (0,1)+(0,k-1)$, but $\aa$ is basic.

\subsubsection{The Alexander dual of the edge ideal of a graph and its fiber cone}
Let $S=\KK[x_1, \ldots ,x_n]$ be the polynomial ring on $n$ variables over a field $\KK$, and let $t$ be an extra indeterminate. Let $G$ be a graph on $n$ vertices. Denote by
\[A(G)_k \subseteq S[t]\]
the infinite dimensional $\KK$-vector space generated by the
\[ \{x_1^{a_1}\cdots x_n^{a_n}t^k  \quad \hbox{ such that $\aa$ is a $k$-cover of } G \}.\]
The graded $S$-algebra
\[A(G)=\bigoplus_{k\geq 0} A(G)_k\] is called the \textit{vertex cover algebra} of $G$.

If $G$ is bipartite and has at least one edge, it admits no indecomposable $k$-covers for $k \geq 2$; thus $A(G)$ is a standard graded $S$-algebra \cite{HHT}. In particular, the quotient
\[\bar{A}(G)=A(G) /\mm A(G)\] is a standard graded $\KK$-algebra. It is easy to check that:

\begin{compactitem}
\item $A(G)$ is generated as a $S$-algebra by the monomials $x_1^{a_1} \cdots x_n^{a_n}t$ such that $\aa=(a_1, \ldots, a_n)$ is an indecomposable $1$-cover;
\item $\AG$  is generated as a $\KK$-algebra by the monomials $x_1^{a_1} \cdots x_n^{a_n}t$ such that $\aa=(a_1, \ldots, a_n)$ is an indecomposable $1$-cover;
\item a basis of $\bar{A}_k(G)$ as a $\KK$-vector space is given by the monomials $x_1^{a_1} \cdots x_n^{a_n}t^k$ such that $\aa=(a_1, \ldots, a_n)$ is a basic $k$-cover;
\item in particular, the Hilbert function $HF_{\AG} (k)=\dim_{\KK} \bar{A}(G)_k$ counts the number of basic $k$-covers of $G$.
\end{compactitem}

Recall that for a homogeneous ideal $I \subseteq S$ the dimension of its fiber cone $F_{\mm}(I)$ is called the \textit{analytic spread} of $I$. When $\KK$ is an infinite field, the analytic spread is the minimal number of generators of any ideal $J$ satisfying
\begin{compactitem}[--]
\item $J$ is contained in $I$;
\item there exists a positive integer $s$ such that $J^s \cdot I = I^{s+1}$;
\item $J$ is minimal among the ideals that satisfy the previous two properties.
\end{compactitem}
In particular, the analytic spread is larger or equal than the arithmetical rank, defined as
\[\ara(I) := \min \{ \hbox{ number of generators of $J$ s.t. } \sqrt{J}=\sqrt{I}\}. \]

Now, let $G$ be a bipartite graph with at least one edge; the ideal
\[I=\bigcap_{ \{i,j\} \hbox{ edge of } G} \left( x_i, \; x_j \right) \; ;\]
is a squarefree monomial ideal, whose fiber cone coincides with $\bar{A}(G)$. In other words, $\bar{A}(G)$ is the fiber cone of the Alexander dual of the edge ideal of $G$, and $\dim (\bar{A}(G))$ is its analytic spread. (Also, $A(G)$ corresponds to the Rees algebra $S[I \,t]$.)

\par
\bigskip
\noindent \textbf{Notes.}
\begin{compactenum}
\item If $G$ is not bipartite, it is still possible to define $\bar{A}(G)$ as the symbolic fiber cone of $I$, provided $G$ has at least one edge; see the last section for a broader approach involving hypergraphs.
\item If $G$ is the disjoint union of $n$ points, $G$ is bipartite but $\AG$ is not finitely generated as  $\KK$-algebra. The basic $k$-covers of $G$ are those that assign $1$ to a single vertex and $0$ to all other vertices; twice such a basic cover is not a basic cover anymore.
\item \textbf{Throughout this paper we always assume that G has at least one edge}.
\end{compactenum}


\section{When is $\bar{A}(G)$ a domain?}
In this paragraph, we answer the following questions: When does $\bar{A}(G)$ contain
\begin{compactitem}[--]
\item zero-divisors?
\item non-zero-divisors?
\item nilpotent elements?
\end{compactitem}

The idea, as in Herzog--Hibi--Trung\cite{HHT},  is to use vertex covers as a ``bridge'' between commutative algebra and combinatorics:
\begin{compactitem}
\item first we characterize the algebraic properties above with statements in the language of vertex covers;
\item then we translate the vertex covers statements into combinatorial properties of the given graph $G$.
\end{compactitem}
This way we establish that $\bar{A}(G)$ is a domain if and only if the graph $G$ has at least one edge and satisfies the following, purely combinatorial condition.

\begin{defn} A bipartite graph $G$ satisfies the \textit{weak square condition} (shortly, \textit{WSC}) if for every non-isolated vertex $i$ of $G$ there exists an edge $\{i,j\}$ of $G$ such that the following holds true:
\[ \hbox{ If } \{i,i'\} \hbox{ is an edge of $G$ and } \{j,j'\} \hbox{ is an edge of $G$, then } \{i',j'\} \hbox{ is also an edge of $G$.}\]
\end{defn}
It is easy to check that any complete bipartite graph satisfies WSC; furthermore, if two graphs satisfy WSC, their disjoint union satisfies WSC as well.

We will also prove the following prospect:
\begin{center}
\begin{tabular}{l|l|l}
$\bar{A}(G)$ has a non-zero divisor & if and only if &$G$ has at least one edge;\\
$\bar{A}(G)$ has a zero-divisor & if and only if & $G$ does not satisy WSC;\\
$\bar{A}(G)$ has a nilpotent element & if and only if &$G$ has no edges.
\end{tabular}
\end{center}

The second answer is particularly interesting in connection with the following result.
\begin{prop} \label{thm:graziehochster}
Suppose that $\bar{A}(G)$ is a domain. Then $\bar{A}(G)$ is a normal domain and a Cohen-Macaulay algebra, too.
\end{prop}

\begin{proof}
If $\AG$ is a domain then it is a semigroup ring. More precisely $\AG \cong \KK[C]$ where $C \subseteq \ZZ^{n+1}$ is the semigroup
\[C=\{(a_1, \ldots ,a_n,k): \ k\geq 1, \  (a_1, \ldots ,a_n) \mbox{ is a basic $k$-cover} \}.\]

By a result of Hochster (\cite[Theorem 1]{Ho}), it suffices to show that $C$ is a normal semigroup. Choose an $\aa \in \ZZ \cdot C$ such that $m \cdot \aa \in C$ for some $m \in \NN$. Write
\[\aa =(\sum_{i=1}^p m_ia_1^i, \ \ldots , \ \sum_{i=1}^p m_ia_n^i, \ \sum_{i=1}^p m_i k_i) \]
where $m_i\in \ZZ$ and $\aa^i=(a_1^i, \ldots, a_n^i, k_i)\in C$ for every $i=1, \ldots,p$. For every edge $\{i,j\}$ of $G$
\begin{center}$m \cdot (\sum_{h=1}^p  m_k(a_i^h+a_j^h)) \ \geq \ m \cdot (\sum_{h=1}^p m_hk_h) \geq 1 \ \implies \
\sum_{h=1}^p  m_h(a_i^h + a_j^h) \ \geq \ \sum_{h=1}^p m_hk_h \geq 1$,
\end{center}
and equality on the left hand side implies equality on the right and side.
So $m \cdot \aa \in C \implies \aa \in C$, then $C$ is a normal semigroup.
\end{proof}

Let us start with some easy lemmas.

\begin{lem} \label{thm:technicalshit}
Let $G$ be a bipartite graph. Then $\bar{A}(G)$ is not a domain if and only if there exist $m$ basic $1$-covers that add up to a non-basic $m$-cover, for some $m \geq 2$.
\end{lem}
\begin{proof}
The {\it if}-part is easy.\\
Conversely, choose two non-zero elements $f,g \in \bar{A}(G)$ such that $f \cdot g =0$; or, in other words, choose two elements $F,G \in A(G) \setminus \mm A(G)$ such that $F \cdot G = H$ with $H \in \mm A(G)$. Consider a term order $\prec$ on $S[t]$. We can assume that the leading term of $F$, $LT_{\prec}(F)$, does not belong in $\mm A(G)$: otherwise, replace $F$ with $F'=F-\alpha \cdot LT_{\prec}F$ for a suitable $\alpha \in \KK$. Similarly, we can also assume $LT_{\prec}(G) \notin \mm A(G)$. On the other hand $\mm A(G)$ is generated as a $\KK$-vector space by monomials, so $LT_{\prec}(H)\in \mm A(G)$.\\
Then we must have $LT_{\prec}(F)\cdot LT_{\prec}(G)=LT_{\prec}(H)$, so
\[ (x_1^{a_1}\cdots x_n^{a_n} \cdot t^k)\cdot (x_1^{b_1}\cdots x_n^{b_n} \cdot t^h)=x_1^{c_1}\cdots x_n^{c_n} \cdot t^m,\]
where $\aa$ is a basic $k$-cover, $\bb$ is a basic $h$-cover, $\cc$ is a non basic $m$-cover,  and $m=k+h$. Find $k$ basic $1$-covers that add up to $\aa$, and $h$ basic $1$-covers that add up to $\bb$. In total we have $m$ basic $1$-covers adding up to some non-basic $m$-cover $\cc$.
\end{proof}

\begin{lem} \label{thm:facile}
If $\bb$ is a basic $k$-cover of a graph $G$ on $[n]$, $b_i \leq k$ for all $i \in \{ 1, \ldots, n\}$.
\end{lem}

\begin{proof}
Whenever $b_i \geq k+1$, $\bb$ can be lopped at $i$.
\end{proof}

\begin{prop} \label{thm:isolated}
Let $G$ be a bipartite graph. Then, $\bar{A}(G)$ has no nilpotent element.
\end{prop}

\begin{proof} If  $\aa$ is a basic $k$-cover of $G$ and $t$ is any natural number, then $t \cdot \aa$ is a basic $(tk)$-cover of $G$ (by our assumption $G$ has at least one edge). In particular, $\bar{A}(G)$ has no nilpotent generator. Now arguing as in Lemma \ref{thm:technicalshit} we deduce that there are no nilpotent elements in $\AG$.
\end{proof}

\begin{os} \label{thm:noisolated}
Let $G$ be a bipartite graph. Suppose $G$ has isolated vertices. Any basic cover of $G$ assigns $0$ to all the isolated vertices; therefore, $\bar{A}(G)$ is isomorphic as a graded  $\KK$-algebra to $\bar{A}(G')$ - where $G'$ is the graph obtained deleting from $G$ all the isolated vertices. In particular $\AG$ is a domain if and only if $\bar{A}(G')$ is a domain.
\end{os}

\begin{thm} \label{thm:graph2covers}
Let $G$ be a bipartite graph. Then, $\bar{A}(G)$ is a domain if and only if for every non-isolated vertex $i$ of $G$ there exists an edge $\{i,j\}$ of $G$ such that for every basic $1$-cover $\cc$ one has $c_i+c_j=1$.
\end{thm}

\begin{proof}
``$\Leftarrow$''. By contradiction, if $\bar{A}(G)$ is not a domain, there exist $m$ basic 1-covers $\aa^1$, $\ldots$,~$\aa^m$ that add up to a non-basic $m$-cover (because of Lemma \ref{thm:technicalshit}). Fix a vertex $i$ of $G$; by the assumption, there is an edge $\{i,j\}$ such that $a^h_i + a^h_j=1$, for all $h \in \{1, \ldots, m\}$. But then the $m$-cover $\aa:=\aa^1 + \ldots + \aa^m$ cannot be lopped at $i$. Repeat this reasoning for all $i$ to conclude that $\aa$ is basic, a contradiction.

``$\Rightarrow$''. We can assume by the previous remark that $G$ has no isolated points. By contradiction, there is a vertex $i$ such that for each vertex $j$ adjacent to $i$, there is some basic $1$-cover $\bb=\bb(j)$  such that $b(j)_i + b(j)_j = 2$. (Notice that $b(j)_i + b(j)_j$ cannot be $3$ or more, because of Lemma \ref{thm:facile}.)
Now define
\[\bb : = \sum_{j \hbox{ adjacent to } i} \bb(j).\]
Clearly $\bb$ is a $d$-cover, where $d>0$ is the degree of the vertex $i$ in the graph $G$. It is also the sum of $d$ basic $1$-covers. It is straightforward to check that $\bb$ can be lopped at $i$, which ends the proof.
\end{proof}

\begin{lem} \label{thm:edge}
Let $\{i,j\}$ be an edge of a bipartite graph $G$. The following are equivalent:
\begin{compactenum}[(I).]
\item $a_i + a_j = 1$ for every basic $1$-cover $\aa$ of $G$;
\item if $\{i, i'\}$ is an edge and $\{j, j'\}$ is an edge, then $\{i', j'\}$ is an edge.
\end{compactenum}
\end{lem}

\begin{proof}
$(I) \Rightarrow (II)$: By contradiction, suppose $\{i',i\}$,$\{i,j\}$, $\{j,j'\}$ are edges of $G$, while $\{i', j'\}$ is not an edge. We can define a $1$-cover $\bb$ as
\[b_h := \left \{ \begin{array}{ll}
                 0 & \hbox{ if } h \in \{ i', j'\} , \\
        1& \hbox{ otherwise.}
\end{array} \right.\]
In general $\bb$ is not basic; however, we can lop it at vertices not in the set $\{i,j,i',j'\}$, until we obtain a basic $1$-cover $\aa$ such that $a_i + a_j = 2$; a contradiction with $(I)$.

$(II) \Rightarrow (I)$: By the assumption, there is an edge $\{i,j\}$ such that for all $i'$ adjacent to $i$, and for all $j'$ adjacent to $j$, $\{i', j'\}$ is also an edge of $G$. By contradiction, suppose there exists a basic $1$-cover $\cc$ of $G$ such that $c_i = c_j =1$. Since $\cc$ is basic, $\cc$ cannot be lopped at $i$, so there must be some vertex $i'$ adjacent to $i$ such that $c_{i'}=0$. Similarly, $\cc$ cannot be lopped at $j$, so there is a $j'$ adjacent to $j$ such that $c_{j'}=0$. But $\{i', j'\}$ is an edge of $G$, and $c_{i'} = c_{j'} =0$: a contradiction.
\end{proof}

\begin{thm}\label{thm:domain}
Let $G$ be a bipartite graph. Then the following are equivalent:
\begin{compactenum}
 \item $G$ satisfies $WSC$;
\item for every non-isolated vertex $i$ of $G$ there exists an edge $\{i,j\}$ of $G$, such that for every basic $1$-cover $\cc$ one has $c_i+c_j=1$;
\item $\bar{A}(G)$ is a domain.
\end{compactenum}
\end{thm}

\begin{proof}
By Remark \ref{thm:noisolated}, it is sufficient to show the equivalence for graphs $G$ without isolated vertices. The equivalence of 2 and 3 is given by Theorem \ref{thm:graph2covers}; the equivalence of 1 and 2 follows directly by Lemma \ref{thm:edge} and the definition of WSC.
\end{proof}

\begin{cor}
 If $G$ satisfies WSC, $\bar{A}(G)$ is a normal Cohen-Macaulay domain.
\end{cor}

\begin{ess}
\begin{compactitem}
\item Let $G$ be a square. Then,  $\bar{A}(G)$ is a domain. If $G$ is a $2k$-cycle, and $k \geq 3$, $G$ is not a domain.
\item Let $G$ be (the 1-skeleton of) a cube. Then,  $\bar{A}(G)$ is not a domain. However, if $\tilde{G}$ is the cube together with its three main diagonals, $\bar{A}(\tilde{G})$ is a domain.
\item If $G$ is a path of $3$ consecutive segments (i.e. a square minus an edge),  $\bar{A}(G)$ is a domain; if $G$ is a path of $4$ consecutive segments (i.e. a pentagon minus an edge),  $\bar{A}(G)$ is not a domain.
\end{compactitem}
\end{ess}

Recall that $G$ is \textit{unmixed} if all the basic $1$-covers of $G$ have the same number of ones. It follows from the definition that when $G$ is unmixed, $m$ basic $1$-covers add up always to a basic $m$-cover: so unmixed implies domain by Lemma \ref{thm:technicalshit}. The converse does not hold: the algebras of complete bipartite graphs are all domains, while $K_{a, b}$ is unmixed if and only if $a=b$.

\begin{prop} \label{thm:domaintrees}
Let $T$ be a tree on $n$ vertices.
\begin{compactitem}
\item $T$ is a domain if and only if every vertex of $T$ is either a leaf, or adjacent to a leaf.
\item $T$ is unmixed if and only if $n$ is even, and $T$ is obtained from some tree $U$ on $\frac{n}{2}$ vertices by attaching one leaf to each vertex of $U$.
\end{compactitem}
\end{prop}

The second part of Proposition \ref{thm:domaintrees} is known: for example, it follows from the result of Herzog--Hibi \cite[Theorem 3.4]{HH} and from the main theorem by Herzog--Hibi--Zheng \cite{HHZ}. Below, however, we give a purely combinatorial proof:

\begin{proof}
The first item follows applying the WSC property to paths: fix a vertex $i$ and consider any edge $\{i,j\}$. Unless $i$ is a leaf (or $j$ is a leaf), we can find edges $\{i,i'\}$ and $\{j,j'\}$: and since a tree is acyclic,  $\{i', j'\}$ is not an edge.

Suppose $T$ is unmixed. Then, the bipartition of the vertices of $B$ gives rise to two basic $1$-covers with the same number of ones: so $n$ is even and $|A|=|B|=\frac{n}{2}$. Moreover, let $\cc$ be the $1$-cover obtained assigning $0$ to all leaves and $1$ to all other vertices. A priori $\cc$ might be not basic; however, since unmixed implies domain, by the previous item every vertex of $T$ is either a leaf or adjacent to a leaf; this implies that $\cc$ is basic.
Let $C$ be the set of elements that are adjacent to at least one leaf. Clearly $\cc$ is identically $1$ on
$C$. So the number of elements of $C$ is $\frac{n}{2}$ (by definition of unmixed). Now we notice that every leaf is adjacent to exactly one element in $C$, while an element $v$ in $C$ is adjacent to $d_v$ leaves; but since $|C|=\frac{n}{2}=n - \#\{ \hbox{ leaves }\}$, by the pidgeon-hole principle we get that $d_v=1$ for all $v$'s.

Conversely, let $T$ be a tree obtained from $U$ as described above, and let $\bb$ be any basic $1$-cover of $U$. $\bb$ can be extended in an unique way to a basic $1$-cover $\bb'$ of $T$ - and all basic $1$-covers of $T$ are obtained this way. The extension consists in assigning $0$ to a leaf that starts at a point that has been given $1$, and the other way around. This means that the obtained basic $1$-cover of $T$ will have exactly $\frac{n}{2}$ ones.
\end{proof}

\begin{es}
Let $T$ be the tree on six vertices and the following five edges:
\[ \{1,2\}, \; \{1,3\}, \; \{1,4\}, \; \{4,5\}, \; \{4,6\}.\]
Then, $T$ is not unmixed, albeit $T$ is a domain.\\
Notice that in this case $I=\cap_{\{i,j\}}(x_i,x_j)=$ $(x_1x_4,$ \ $x_1x_5x_6,$ \ $x_2x_3x_4,$ \ $x_2x_3x_5x_6)$. So $I$ is not generated in only one grade neither a complete intersection, albeit its fiber cone is a domain.
\end{es}

We finally give a criterion to determine whether a basic $k$-cover is a zero-divisor in the algebra or not.

\begin{lem} \label{thm:specialnonzero}
Let $G$ be a bipartite graph, and $\bb$ a basic $k$-cover. If $b_i + b_j = k$ for each edge $\{i,j\}$ of $G$, then the sum of $\bb$ and other basic covers gives always a basic cover.
\end{lem}

\begin{proof}
Let $\aa = \bb+\cc$ be the sum of the basic $k$-cover $\bb$ above and an arbitrary basic $k'$-cover $\cc$. $\aa$ is a $(k+k')$-cover; we claim that $\aa$ cannot be lopped at any vertex. By contradiction, suppose $\aa$ can be lopped at $i$: then, the inequality $a_j +a_i \geq k+k'$ was strict for all $j$'s adjacent to $i$. Besides, $b_i + b_j = k$ for all $j$'s adjacent to $i$. So, $c_j +c_i \geq k'$ must be strict for all $j$'s adjacent to $i$ - which means that $\cc$ can be lopped at $i$, a contradiction.
\end{proof}

\begin{prop}\label{thm:depth}
Let $G$ be a bipartite graph with at least one edge. Then $\depth \left( \AG \right) \geq 2$.
\end{prop}

\begin{proof}
By Remark \ref{thm:noisolated} we  can  assume that $G$ has no isolated vertices. Under this assumption, the vectors $\aa$ and $\bb$ defined by
\[a_h := \left \{ \begin{array}{ll}
                 1 & \hbox{ if } h \in A, \\
        0& \hbox{ if } h \in B
\end{array} \right.\]
\[b_h := \left \{ \begin{array}{ll}
                 0 & \hbox{ if } h \in A, \\
        1& \hbox{ if } h \in B
\end{array} \right.\]
are both basic $k$-covers, where $[n]=A \cup B$ is the bipartition of the vertex set of $G$.

We claim that $\mathbf{x}^{\aa}t=x_1^{a_1} \cdots x_n^{a_n}t$ and $\mathbf{x}^{\bb}t=x_1^{b_1} \cdots x_n^{b_n}t$ form a regular sequence in the algebra $\bar{A}(G)$.

In fact, by Lemma \ref{thm:specialnonzero}, the sum of $\aa$ and other basic covers gives always a basic cover, so $\mathbf{x}^{\aa}t$ is a non-zero divisor arguing like in the proof of Lemma \ref{thm:technicalshit}. To prove that $\mathbf{x}^{\bb}t$ is not a zero-divisor of $\AG/\mathbf{x}^{\aa}t \cdot \AG$ we assume $\mathbf{x}^{\bb}t \cdot f \in (x_1, \ldots ,x_n, \mathbf{x}^{\aa}t)\AG$ and we show $f \in (x_1, \ldots ,x_n, \mathbf{x}^{\aa}t)\AG$ (here $f \in A(G)$ and the operations are in $S$). Again, arguing as in Lemma \ref{thm:technicalshit} we can assume that $f=x_1^{c_1}\cdots x_n^{c_n}t^k$ where $\cc=(c_1, \ldots ,c_n)$ is a basic $k$-cover. Moreover, since $\mathbf{x}^{\bb}t$ is not a zero-divisor (again by Lemma \ref{thm:specialnonzero}) we have that there exists a basic $k$-cover $\dd$ such that
$\aa + \dd=\bb+\cc$. Now it is easy to see that, by the choice of $\aa$ and $\bb$, this implies that $\cc=\aa+\cc'$ for a suitable basic $(k-1)$-cover $\cc'$, and this allows us to conclude.
\end{proof}

\begin{cor}
Let $G$ be a bipartite graph.
Then the projective scheme $\Proj(\AG)$ is connected.
\end{cor}

\begin{proof}
From a result of Hartshorne (see \cite[Proposition 2.1]{Ha}) one has that for any projective scheme $X=\Proj(R)$ the following fact holds: $\depth(R)\geq 2 \implies X$ connected.
\end{proof}

\begin{thm}[Characterization of non-zero-divisors.]
Let $G$ be a bipartite graph, and $\bb$ a basic $k$-cover.
The following are equivalent:
\begin{compactenum}
 \item the sum of $\bb$ and other basic covers gives always a basic cover;
\item for each edge $\{i,j\}$ of $G$ such that $b_i + b_j > k$, there exist an $i'$ adjacent to $i$ and a $j'$ adjacent to $j$ such that for each basic 1-cover $\aa$, one has $a_i + a_{i'} = 1 = a_j + a_{j'}$;
\item for each edge $\{i,j\}$ of $G$ such that $b_i + b_j > k$, there exist an $i'$ adjacent to $i$ and a $j'$ adjacent to $j$ such that:
    \begin{compactitem}[--]
     \item if $\{v, i\}$ is an edge and $\{v', i'\}$ is an edge, then $\{v, v'\}$ is also an edge;
    \item if $\{w, j\}$ is an edge and $\{w', j'\}$ is an edge, then $\{w, w'\}$ is also an edge.
    \end{compactitem}
 \end{compactenum}
\end{thm}

\begin{proof}
The equivalence of 2 and 3 is easily shown via Lemma \ref{thm:edge}.

$1 \Rightarrow 2$: By contradiction, let $\{i, j\}$ be an edge with $b_i + b_j > k$, and suppose that for all $i'$ adjacent to $i$ and $j'$ adjacent to $j$, one can find either a basic $1$-cover $\cc=\cc(i')$ (depending on $i'$) such that $c_i + c_{i'}=2$, or a basic 1-cover $\dd=\dd(j')$ (depending on $j'$)  such that $d_j + d_{j'}=2$.

Consider
\[ \aa:= \sum_{i' \hbox{ adjacent to } i} \cc(i') \; + \; \sum_{j' \hbox{ adjacent to } j} \dd(j'). \]

This $\aa$ is a $h$-cover, where $h$ is the sum of the degrees of $i$ and $j$ in the graph $G$. It is easy to see that $\aa$ can be lopped either at $i$ or at $j$.

$2 \Rightarrow 1$: By contradiction, let $\cc$ be a basic $k'$-cover such that $\aa=\bb+\cc$ is not basic. Reasoning as in Lemma \ref{thm:specialnonzero}, it is easy to see that $\cc$ could be lopped only at those vertices $i$ such that $b_i + b_j > k$ for some $j$ adjacent to $i$. So, fix such an edge $\{i, j\}$ with $b_i + b_j > k$, so that $\cc$ can be lopped at $i$. Clearly $\cc$ is a sum of $k'$ indecomposable $1$-covers; by the assumption, then, $c_i + c_{i'} = k' = c_j + c_{j'}$, for some $i'$ adjacent to $i$, and some $j'$ adjacent to $j$. But if $\cc$ can be lopped at $i$, the inequality $c_i + c_{i'} \geq k'$ should be strict: a contradiction.
\end{proof}

\section{The dimension of $\bar{A}(G)$}

In this paragraph we study the dimension of $\AG$, which is the analytic spread of the ideal
\[I=\bigcap_{\{i,j\} \mbox{ edge of }G}(x_i,x_j)\subseteq S=\KK[x_1, \ldots ,x_n].\]

As we have already said the Hilbert function of $\AG$
\[HF_{\bar{A}(G)}(k)=\dim_{\KK} \bar{A}(G)_k\]
counts the number of basic $k$-covers of $G$. For $k$ large, the Hilbert function grows as $C \cdot k^{d}$, where $d=\dim \AG -1$ and $C\cdot d!=e(\AG)$, the multiplicity of $\AG$. Hence the main tool for computing the dimension will be approximating the number of basic $k$-cover of the graph for $k$ large.

We introduce herein the ``graphical dimension'', an integer that depends only on the combinatorics of $G$.
This invariant captured our interest for the following reasons:

\begin{compactitem}[$\qquad \rhd$]
\item $\gdim (G)  \leq \dim \AG$ (Theorem \ref{thm:inequalities});
\item $\gdim$ can be computed or estimated from below; \footnote{Computing $\gdim$ takes ``too long'' if $G$ has many edges. However, any drawing that represents $G$ as bipartite graph basically gives a lower bound for $\gdim(G)$, and thus for $\dim \AG$: See Definition \ref{thm:gdim}.}
\item $\gdim(C_{2a}) = \dim \bar{A}(C_{2a}) = a$ for each $(2a)$-gon (see Proposition \ref{thm:dimcycle});
\item $\gdim(G)=\dim \bar{A}(G) = \lfloor\frac{n}{2} \rfloor + 1$ if $G$ is a path on $n$ vertices, or if $G$ has a Cohen-Macaulay edge ideal (see Corollaries \ref{thm:dimpath} and \ref{thm:dimCM});
\item $\gdim (T) = \dim \bar{A}(T)$ for each tree $T$ (Theorem \ref{6}).
\end{compactitem}

\vspace{1mm}

\noindent Notice that by Remark \ref{thm:noisolated}, we may assume that $G$ has no isolated points. We can even reduce to the case where $G$ is connected by means of the following Lemma:

\begin{lem}\label{connected}
If $G$ is the disjoint union of two graphs $G_1$ and $G_2$,
\[\dim( \bar{A}(G) ) = \dim( \bar{A}(G_1) ) + \dim( \bar{A}(G_2) ) - 1.\]
\end{lem}

\begin{proof}
Every basic $k$-cover of $G$, restricted to $G_i$, gives a basic $k$-cover of $G_i$; and all basic $k$-covers of $G$ arise this way. This means that \[HF_{\bar{A}(G)}(k) \; = \; HF_{\bar{A}(G_1)}(k) \; \; \cdot \; \; HF_{\bar{A}(G_2)}(k),\]
whence for large $k$ we obtain that
\[\dim( \bar{A}(G) ) - 1 \; = \; \dim( \bar{A}(G_1) ) - 1 \;   + \;  \dim( \bar{A}(G_2) ) - 1.\]
\end{proof}

\begin{defn} \label{thm:gdim}
Let $G$ be a (connected) bipartite graph on $n$ vertices, $n \geq 2$. Let $[n]=A \cup B$ be the vertex partition (unique because of connectedness); we assume from now on $a:=|A| \leq |B|=:b$. A \textit{standard drawing} $\Gamma_G$ of $G$ is a drawing of
$G$ so that:
\begin{compactitem}[--]
\item the elements of $A$ are drawn on a horizontal line, and labeled $1, \ldots, a$ from left to right;
\item the elements of $B$ are drawn on a horizontal line below the previous one, and labeled $1', 2', \ldots, b'$;
\end{compactitem}
Any standard drawing of $G$ partitions the edges of $G$ into three categories:
\begin{compactdesc}
\item{$\qquad \bullet \;$ \textit{slash edges,}} or edges from $i \in A$ to $j' \in B$ with $i > j$.
\item{$\qquad \bullet \;$ \textit{backslash edges,}} namely edges from $i \in A$ to $j' \in B$ with $i < j$.
\item{$\qquad \bullet \;$ \textit{vertical edges,}} namely edges from $i \in A$ to $j' \in B$ with $i = j$.
\end{compactdesc}
\end{defn}

\begin{os}
There is an abuse of notation in ``relabeling'' the vertices $\{1, \ldots, n\}$ as $\{1, \ldots, a\}$ and $\{1', \ldots, b'\}$. We apologize with the reader, although we are convinced that once a drawing has been fixed, this abuse does not represent a problem. In Subsection \ref{tree}, however, we let the drawings vary -- which may cause some confusion in the labeling. Therefore, in Subsection \ref{tree} we will introduce of a more precise (yet ``heavier'') notation.
\end{os}

A graph $G$ might have many distinguished standard drawings -- accordingly with how we label the elements in $A$ and in $B$. To each standard drawing $\Gamma_G$ of $G$ we associate an integer $r=r \left(\Gamma_G \right)$  in $\{0, 1, \ldots, a\}$, as follows:
\begin{defn} Fix a drawing $\Gamma_G$ of $G$. Define $r=r(\Gamma_G)$ as the largest integer that satisfies the following two properties:
\begin{compactenum}
\item $\{1, 1'\}, \{2, 2'\}, \ldots, \{r, r'\}$ are vertical edges of $G$;
\item If $i \leq r$, there is no slash edge between $i$ and $j'$, if $j$ is smaller than $i$.
\end{compactenum}
\end{defn}
Different standard drawings give rise to different $r \left(\Gamma_G \right)$'s; we define the \textit{graphical dimension} of $G$ as
\[\gdim(G):= 1 + \max \left\{r \left(\Gamma_G \right) \; : \; \Gamma_G \hbox{ is a standard drawing of } G \right \}.\]
We will say that $\Gamma_G$ is \textit{optimal} if $r(\Gamma_G)+1=\gdim(G)$. From now on we will often write just $\Gamma$ instead of $\Gamma_G$.

The definition of graphical dimension extends to non-connected bipartite graphs consistently with Lemma \ref{connected}:

\begin{defn}
Let $G$ be a bipartite graph (possibly with isolated points).  Let
\[G \; = \; G_1 \; \cup \ldots \cup \; G_m \; \cup \; H_1\; \cup \ldots \cup \; H_l \]
be the decomposition of $G$ into its connected components, so that every $H_i$ consists of a single isolated vertex, while every $G_i$ has at least one edge. (If $G$ has no isolated vertices $l=0$; by our convention $m \geq 1$). Then
 \[ \gdim(G):= 1 - m + \sum_{i=1}^m \gdim(G_i).\]
\end{defn}

Notice that passing from a drawing $\Gamma$ to another drawing $\Gamma'$ we have to apply a permutation $\sigma$ (resp. a permutation $\tau$) on the set $A$ (resp. on the set $B$). We will use this notation in some proof.

\begin{es}\label{alex}
We show that there always exists a connected bipartite graph $G$ for any triple of numbers $a$, $b$, $r=\gdim(G)-1$ such  that $1\leq r \leq a \leq b$.\\
If $r=1$ it is easy to see that $G=K_{a,b}$ is one such graph. So given natural numbers $2\le r \le a \le b$ consider the following graph $G$:

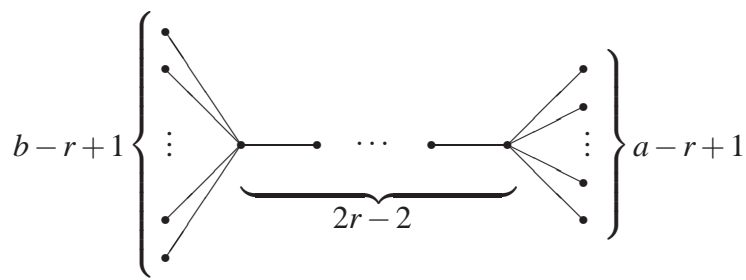
\[
\setlength{\unitlength}{5mm}
\begin{picture}(19,11)
\put(6,6){\circle*{0.2}}
\put(8,6){\circle*{0.2}}
\put(11,6){\circle*{0.2}}
\put(13,6){\circle*{0.2}}
\put(4,3){\circle*{0.2}}
\put(4,4){\circle*{0.2}}
\put(4,8){\circle*{0.2}}
\put(4,9){\circle*{0.2}}
\put(15,4){\circle*{0.2}}
\put(15,5){\circle*{0.2}}
\put(15,7){\circle*{0.2}}
\put(15,8){\circle*{0.2}}

\put(6,6){\line(1,0){2}}
\put(11,6){\line(1,0){2}}
\put(6,6){\line(-2,-3){2}}
\put(6,6){\line(-1,-1){2}}
\put(6,6){\line(-1,1){2}}
\put(6,6){\line(-2,3){2}}
\put(13,6){\line(1,-1){2}}
\put(13,6){\line(2,-1){2}}
\put(13,6){\line(2,1){2}}
\put(13,6){\line(1,1){2}}

\put(9,6){$\ldots$}
\put(4,5.7){$\vdots$}
\put(15,5.7){$\vdots$}
\put(3,5.8){$\left\{ \begin{array}{ll}&\\&\\&\\&\\&\\&\\&\\\end{array}\right.$}
\put(14,5.8){$\left. \begin{array}{ll}&\\&\\&\\&\\&\\\end{array}\right\}$}
\put(6,5){$\underbrace{\phantom{zzzzzzzzzzzzzzzzzzzzzz}}  $}
\put(0,5.7){$b-r+1$}
\put(8.4,3.8){$2r-2$}
\put(16.3,5.7){$a-r+1$}

\end{picture}
\]

As one could not have any more vertical edges, it is easy to see that the following is an optimal standard drawing of $G$ with exactly $a$ vertices above and $b$ vertices below.
.
\[
\setlength{\unitlength}{5mm}
\begin{picture}(17,6)

\put(1,1){\circle*{0.2}}
\put(3,1){\circle*{0.2}}
\put(5,1){\circle*{0.2}}
\put(8,1){\circle*{0.2}}
\put(10,1){\circle*{0.2}}
\put(13,1){\circle*{0.2}}
\put(17,1){\circle*{0.2}}

\put(1,5){\circle*{0.2}}
\put(3,5){\circle*{0.2}}
\put(5,5){\circle*{0.2}}
\put(8,5){\circle*{0.2}}
\put(10,5){\circle*{0.2}}
\put(13,5){\circle*{0.2}}
\put(16,5){\circle*{0.2}}

\put(1,1){\line(0,1){4}}
\put(3,1){\line(0,1){4}}
\put(5,1){\line(0,1){4}}
\put(8,1){\line(0,1){4}}
\put(10,1){\line(0,1){4}}

\put(1,5){\line(1,-2){2}}
\put(3,5){\line(1,-2){2}}
\put(5,5){\line(1,-2){0.3}}
\put(8,1){\line(-1,2){0.3}}
\put(8,5){\line(1,-2){2}}

\put(10,1){\line(3,4){3}}
\put(10,1){\line(3,2){6}}

\put(1,5){\line(3,-1){4.3}}
\put(1,5){\line(4,-1){4.3}}

\put(13,1){\line(-3,1){5.3}}
\put(17,1){\line(-4,1){9.3}}

\put(6,3){$\ldots$}
\put(14,1){$\ldots$}
\put(14,5){$\ldots$}
\put(0.8,0.8){$\underbrace{\phantom{zzzzzzzzzzizzzzzzzzzzzzzzzzzz}}  $}
\put(5.3,-0.3){$r$}
\end{picture}
\]
 Thus, we have that $\gdim(G) = r+1$.
\end{es}

\vspace{2mm}

\begin{es}\label{9}
Let $G=C_{10}$. Here comes a standard drawing for $G$:

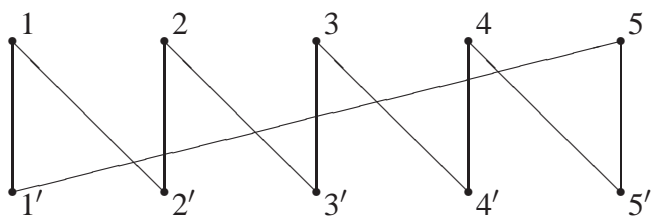
\setlength{\unitlength}{1cm}
\begin{picture}(15,4)(-1.5,-1)
\put(2,0){\circle*{0.1}}
\put(4,0){\circle*{0.1}}
\put(6,0){\circle*{0.1}}
\put(8,0){\circle*{0.1}}
\put(10,0){\circle*{0.1}}

\put(2,2){\circle*{0.1}}
\put(4,2){\circle*{0.1}}
\put(6,2){\circle*{0.1}}
\put(8,2){\circle*{0.1}}
\put(10,2){\circle*{0.1}}

\put(2,0){\line(0,1){2}}
\put(2,0){\line(4,1){8}}
\put(4,0){\line(-1,1){2}}
\put(4,0){\line(0,1){2}}
\put(6,0){\line(-1,1){2}}
\put(6,0){\line(0,1){2}}
\put(8,0){\line(-1,1){2}}
\put(8,0){\line(0,1){2}}
\put(10,0){\line(-1,1){2}}
\put(10,0){\line(0,1){2}}
\put(2.1,2.1){$1$}
\put(4.1,2.1){$2$}
\put(6.1,2.1){$3$}
\put(8.1,2.1){$4$}
\put(10.1,2.1){$5$}
\put(2.1,-0.3){$1'$}
\put(4.1,-0.3){$2'$}
\put(6.1,-0.3){$3'$}
\put(8.1,-0.3){$4'$}
\put(10.1,-0.3){$5'$}

\end{picture}

In the drawing above $r(\Gamma)$ is four, so the $\gdim(G)$ is at least five. It is easy to see that $\gdim(G)$ cannot be six (for example because $G$ has no leaves). Thus $\gdim(G)=5.$
\end{es}

The previous example generalizes as follows.

\begin{prop} \label{thm:drawingcycles}
If $G$ is a $2a$-gon, then the graphical dimension of $G$ is $a$.
\end{prop}

We are now ready to prove that $\gdim$ is smaller or equal than $\dim \AG$.

\begin{lem}\label{thm:3}
Let $\Gamma$ be an optimal standard drawing of $G$. If $i \geq \gdim(G)$, the vertex $i$ in $A$ must be connected with a slash edge to some $h' \in B$ such that $h < \gdim(G)$ in $B$, and the vertex $i'$ in $B$ must be connected with a backslash edge to some $j < \gdim(G)$ in $A$.
\end{lem}

\begin{proof}
We argue by contradiction: were $i$ adjacent only to vertices $j'$ with $j \geq \gdim(G)$, one would get a better drawing as follows:
\begin{compactitem}
\item in $A$, move the vertex $i$ into the $(r(\Gamma) +1)$-th position (from the left to the right):
\item in $B$, move the vertex $j'$ into the $(r(\Gamma) +1)$-th position.
\end{compactitem}
In the new drawing $r$ has increased by one, a contradiction with $\Gamma$ being optimal.

Analogously, if $i'$ in $B$ were only connected with vertices $j \geq \gdim(G)$, one could improve the drawing by moving $j$ into the leftmost position in $A$, and $i'$ into the leftmost position in $B$.
\end{proof}

In view of the previous lemma, for an optimal standard drawing $\Gamma$ of $G$ we can define (set $r=r(\Gamma)$)
\[ m(i, \Gamma)=\min \{j=1, \ldots ,r: \{i,j'\} \mbox{ is an edge}\} \ \ \ \forall \ i=r+1, \ldots, a\]
\[M(i, \Gamma)=\max \{j=1, \ldots ,r: \{i,j'\} \mbox{ is an edge}\} \ \ \ \forall \ i=r+1, \ldots, a\]
and analogously $m(j',\Gamma)$ and $M(j',\Gamma)$ for every $j=r+1, \ldots ,b$.

\begin{thm} \label{thm:inequalities}
Let $G$ be a connected bipartite graph. We have the following inequalities:
\[\gdim(G) \le \dim{\bar{A}(G)} \le a+1\]
\end{thm}

\begin{proof}
We denote $r:= \gdim(G) -1$. We argue by showing

\begin{equation}
\!\! \left \{ \! \!
    \begin{array}{c}
    \hbox{weakly descending sequences}  \\
    \hbox{of $r$ numbers in } \{0, 1, \ldots, k\}
    \end{array}
\right \}
\hookrightarrow
\left \{
    \begin{array}{c}
    \hbox{basic $k$-covers} \\
    \hbox{of the graph } G
    \end{array}
\right \}
\hookrightarrow
\left \{
    \begin{array}{c}
    \hbox{words of $a$ letters in} \\
    \hbox{the alphabet } \{0, 1, \ldots, k\}
    \end{array}
\!\!\!\right \} .
\end{equation}
This is enough to conclude: in fact, passing to the cardinalities we get
\[ \binom{k+r}{r} \; \leq \; HF_{\bar{A}(G)}(k) \; \leq \; (k+1)^a,
\]
and a trivial asymptotic analysis (for $k$ large) of such inequality yields
\[r \; \leq \; \dim \bar{A}(G) - 1 \; \leq \; a.\]
We first focus on the second injective map in (1). Take a basic $k$-cover
\[\aa=\left(a_1, a_2, \ldots, a_a, a_{1'}, a_{2'}, \ldots, a_{b'} \right).\]
By Lemma \ref{thm:facile}, each $a_i$ belongs to $\{0, \ldots, k\}$; so the restriction of $\aa$ to the set $A$ is a word of $a$ letters in the alphabet $\{0, \ldots, k\}$.

Moreover,  as $\mathbf{a}$ is a basic $k$-cover, for every $j = 1,\ldots, b$, there exists an index $i = 1,\ldots,a$ with $\{i,j'\}$ is an edge and $a_i + a_{j'} = k$. Moreover for every $h =1, \ldots ,a$ such that $\{h,j'\}$ is an edge, then $a_h+a_{j'}\geq k$. So for every vertex $j' \in B$ we have
\[ a_{j'} = k - \min\{ a_i \hbox{ such that } \{i,j'\} \hbox{ is an edge of } G\}.\]
Thus every basic $k$-cover $\aa$ is uniquely determined by its restriction to the set $A$, which proves the second inclusion in (1).

As far as the first injective map is concerned, take an optimal standard drawing $\Gamma$ of $G$ and a descending sequence $\omega_1 \geq \ldots \geq \omega_r$ (where $r=r(\Gamma)$) in the alphabet $\{0, 1, \ldots, k\}$. Can we produce a 1-cover from that? Recall that according to our convention, the vertices of $A$ are labeled $1, 2, \ldots, r, \ldots, a$, while the vertices of $B$ are labeled $1', 2', \ldots, r', \ldots, b'$. Thus we define:
 \begin{compactdesc}
 \item{$\quad$ \textit{(above, to the left of $r$)}: } $a_i := \omega_i$ for $i = 1,\ldots,r$;
\item{$\quad$ \textit{(below, to the left of $r$)}: } $a_{j'}:= k - a_j$ for $j=1, \ldots, r$;
 \item{$\quad$ \textit{(above, to the right of $r$)}: } $a_i:= k - a_{m(i,\Gamma)'}$ for $i = r+1,\ldots, a$;
 \item{$\quad$ \textit{(below, to the right of $r$)}: } $a_{j'}:= k - \min \{a_i \hbox{ s. t. } i=1, \ldots ,a \hbox{ and } \{i, j'\} \hbox{ is an edge of } G \},$ for $j = r+1,\ldots, b$.
 \end{compactdesc}
The correctness of this definition follows from Lemma \ref{thm:3}. The weak monotonicity assumption is needed because of the ``backslash edges''. For example, suppose we have a backslash edge from $1$ to $3'$. Since we defined $a_1=\omega_1$ and $a_3'=k -\omega_3$, we should check that $a_1 + a_3 \geq k$ in order for $\aa$ to be a $k$-cover; yet this is guaranteed by the assumption $\omega_1 \geq \omega_3$.

It is instead easy to check that what we defined is a basic $k$-cover of $G$.

\end{proof}

%

As we said in the Introduction, we do not know whether $\dim(\AG)=\gdim(G)$ holds true for generic bipartite graphs as well. However, as it is true for trees and cycles (Theorem \ref{6} and Proposition \ref{thm:dimcycle}) and as so far there is no counterexample, we think that this is true in general.
Here are two partial answers.





\begin{cor}\label{thm:dimpath}
If $G$ is a path, \[\gdim(G) = \dim{\bar{A}(G)} = \left \lfloor \frac{n}{2} \right \rfloor + 1. \]
\end{cor}

\begin{proof}
Draw $G$ as a zig-zag path starting at $1'$ (i.e. with edges $\{1',1\},$ $ \{1,2'\},$ $ \{2',2\},$ $ \{2,3'\}$, and so on). If $n=2a$, or $n=2a+1$, the drawing obtained is a standard drawing $\Gamma$ with $r(\Gamma)=a$;
and by Theorem \ref{thm:inequalities} we obtain
\[a+1 \leq \gdim (G) \leq \dim{\bar{A}(G)} \leq a+1.\]
\end{proof}

Recall that a graph $G$ is said to be Cohen-Macaulay if $S/I(G)$ ($I(G)$ denotes the edge ideal of $G$) is Cohen-Macaulay for every field $\KK$.

\begin{cor}\label{thm:dimCM}
Let $G$ be a Cohen-Macaulay bipartite graph with $n$ vertices. Then \[\gdim(G) = \dim{\bar{A}(G)} = \frac{n}{2} + 1. \]
\end{cor}

\begin{proof}
By a result of Herzog--Hibi \cite[Theorem 3.4]{HH} $n$ must be even and $G$ must ``come from a poset $P$ of $\frac{n}{2}$ elements'', meaning that $G$ is obtained by drawing a set $A$ of $\frac{n}{2}$ points, a set $B$ of $\frac{n}{2}$ points, and some edges between them, so that the $i$-th leftmost vertex of $A$ and the $j$-th leftmost vertex of $B$ are connected by an edge if and only if $p_i \prec p_j$ in the poset $P$. Of course we can assume that $P$ is labeled so that $p_i \prec p_j$ implies $i \leq j$. This way, if $G$ is connected, the drawing above is automatically a standard drawing of $G$, with $r(\Gamma)=\frac{n}{2}= a$; but then \[a+1 \leq \gdim (G) \leq \dim{\bar{A}(G)} \leq a+1,\] and we are done.

If $G$ is not connected, the connected components of $G$ are also Cohen-Macaulay, so the statement follows by Lemma \ref{connected}.
\end{proof}

We can strengthen the second inequality of Theorem \ref{thm:inequalities} as follows.

\begin{thm}\label{thm:7}
Let $G$ be a bipartite graph. Suppose that for every $j'\in B$ the degree of the vertex $j'$ in $G$ is not smaller than $s$.  Then
\begin{compactitem}[--]
\item $\dim(\AG)\leq a-s+2$;
\item if equality holds above, $\e(\AG)\leq \displaystyle \binom{a}{s}(a-s)!$.
\end{compactitem}
\end{thm}
\begin{proof} The idea is to show that
\begin{equation}
\left \{
    \begin{array}{c}
    \hbox{basic $k$-covers} \\
    \hbox{of the graph } G
    \end{array}
\right \}
\hookrightarrow
\left \{
    \begin{array}{c}
    \left(\omega_1, \omega_2, \ldots, \omega_a \right) \hbox{ in } \{0, 1, \ldots, k\}^a \\
    \hbox{ with } \# \{ i ~|~  \omega_i \hbox{ is maximal } \} \geq s
    \end{array}
\!\!\!\right \}.
\end{equation}
In fact, let $\aa$ be a basic $k$-cover; denote $M:=\max \{a_i:i\in A\}$ and $A_M:=\{i\in A:a_i=M\}$. By contradiction suppose that $\# A_M < s$. Then by considering an $i \in A_M$, we have that for every $j'\in B$ with $\{i,j'\}\in E$ there is an $i_j$ such that $a_{i_j}<M=a_i$ and $\{i_j,j'\}\in E$ (because the degree of $j'$ is greater or equal to $s$). So $a_{j'}\geq k-a_{i_j}$ because $\aa$ is a $k$-cover. But then $a_i+a_{j'}>M$ for every $j' \in B$ such that $\{i,j'\}\in E$, a contradiction with $\aa$ being basic. So the map in (2) is injective.

Let us compute the cardinality of the set on the right-hand side in (2). First fix $a$ and $s$. How many $a$-tuples of integers in $\{0, \ldots, k\}$ are there, such that their global maximum $M$ is attained by \textit{exactly} $t$ of them? The answer is:
\[ \sum_{M=1}^k \binom{a}{t} (M-1)^{a-t}.\]
In fact, any such $a$-tuple is obtained fixing the maximum $M$, choosing the $t$ positions out of $\{1, 2, \ldots, a\}$ in which this maximum should be attained, and then filling the remaining $a-t$ positions with arbitrary numbers in $\{0, \ldots, M-1\}$.

Yet the cardinality of the right-hand side in (2) is the number of $a$-tuples in $\{0, \ldots, k\}$  such that their global maximum $M$ is attained by \textit{at least} $s$ of them. This is of course counted by
\[ \sum_{t=s}^{a} \; \sum_{M=1}^k \binom{a}{t} (M-1)^{a-t},\]
which after the substitutions $i=M-1$ and $m=a-t$ becomes
\[ \sum_{m=0}^{a-s} \binom{a}{m} \sum_{i=0}^{k-1}  i^{m}.\]
Now, it is well known that for any fixed $m \in \NN$, $\displaystyle \sum_{i=0}^{k-1}  i^{m}$ is a polynomial (the so-called ``Eulerian polynomial'') with leading term \[\frac{1}{m+1} \; k^{m+1}.\] It follows that
\[ HF_{\bar{A}(G)}(k) \; \leq \; \binom{a}{s} \cdot \frac{1}{a-s+1} \cdot k^{a-s+1} + O(k^{a-s}).\]
\end{proof}

\begin{cor}
For any bipartite graph $G$, the following are equivalent:
\begin{compactitem}
\item $\dim(\AG)=2$;
\item $G=K_{a,b}$ for some positive integers $a, \ b$;
\item $\AG$ is a polynomial ring in two variables over $\KK$.
\end{compactitem}
\end{cor}
\begin{proof}
If $\dim(\AG)=2$ by Theorem \ref{thm:inequalities} $\gdim(G)\leq 2$, and it is easy to show that this is possible only if $G=K_{a,b}$. By Theorem \ref{thm:7} and Theorem \ref{thm:inequalities} $\dim(\bar{A}(K_{a,b}))=\gdim(K_{a,b})=2$; moreover $\bar{A}(K_{a,b})\cong P/J$ where $P$ is a polynomial ring over $\KK$ in $N=\#\{\mbox{basic 1-cover of }K_{a,b}\}$ variables and $J \subseteq P$ is a homogeneous ideal (see the last section for an explanation of this). As there are only two basic 1-cover of $K_{a,b}$, $J=0$ and we are done.
\end{proof}

\begin{prop} \label{thm:dimcycle}
Let $G=C_{2a}$ be an even cycle.
Then,
\begin{compactitem}[--]
\item $\gdim(G) = \dim(\AG)=a$;
\item $\e(\AG)\leq \displaystyle \binom{a}{2}(a-2)!$
\end{compactitem}
\end{prop}

\begin{proof}
We know from Proposition \ref{thm:drawingcycles} that $\gdim(G) =a$. By Theorem \ref{thm:inequalities}, then, $a \leq \dim \AG \leq a+1$. Besides, all the vertices of $G$ have degree two: so we can apply Theorem \ref{thm:7} with $s=2$ and obtain that $\dim \AG \leq a$. Therefore $\gdim(G) = \dim(\AG)=a$; and reapplying Theorem \ref{thm:7}, we get that $\e(\AG)\leq \displaystyle \binom{a}{2}(a-2)!$.
\end{proof}

\begin{es}
One can prove that if $G=C_6$, the inclusion of Theorem \ref{thm:7} is actually an equality; hence, the multiplicity of $\bar{A}(C_6)$ is exactly 3.\\
Moreover $\bar{A}(G)$ can be presented as a quotient of the polynomial ring $P=\KK[X_1, \ldots ,X_N]$ where $N$ is the number of basic 1-covers of $G$. So
\[ \bar{A}(G)\cong P/J \]
where $J$ is a homogeneous ideal of $P$ contained in $(X_1,\ldots ,X_N)^2$.
For the hexagon $N=5$, hence $\operatorname{height}(J)=N-\dim(\bar{A}(C_6))=2$. Then the multiplicity of $\bar{A}(C_6)$ is minimal, i.e. $e(\bar{A}(C_6))=\operatorname{height}(J)+1$.\\
\end{es}
The example above can be generalized as follows: recall that a graph is said $d$-regular if every vertex has degree $d$.

\begin{prop} \label{min multiplicity}
Let $a \geq 2$ be an integer. Let $G$ be an $(a-1)$-regular bipartite graph with $n=2a$ vertices. Then $\dim(\AG)=3$ and $\AG \cong P/J$ has minimal multiplicity.
\end{prop}
\begin{proof}
Notice that the basic $1$-covers of $G$ are exactly $a+2$:
\begin{compactitem}
\item the unique basic 1-cover which assigns 0 to all the vertices of $A$;
\item the unique basic 1-cover which assigns 1 to all the vertices of $A$;
\item for every fixed vertex $i$ of $A$, the unique basic 1-cover which assigns 0 to $i$ and 1 to all the other vertices of $A$.
\end{compactitem}
By Theorems \ref{thm:inequalities} and \ref{thm:7} it follows that $\dim(\AG)=3$. So $P$ is a polynomial ring of $a+2$ variables, and $\operatorname{height}(J)=a-1$. Now one can show that in this case the inclusion of Theorem \ref{thm:7} is actually an equality; hence 
\[e(\AG)=a=\operatorname{height}(J)+1.\]
\end{proof}

Now we show an example of a bipartite graph $G$ such that $\AG$ is not Cohen-Macaulay.

\begin{es}\label{G not CM} (Non-Cohen-Macaulay Algebra).
Let $G$ be the path in the picture below (notice that $\gdim(G) =4$):

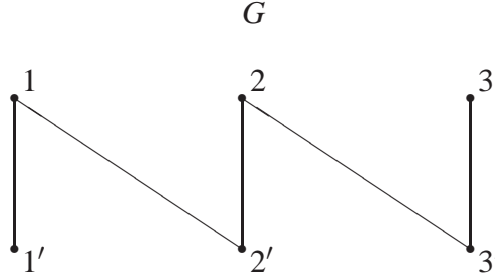
\setlength{\unitlength}{1cm}
\begin{picture}(15,5)(0,-1)
\put(6,3){$G$}

\put(3,0){\circle*{0.1}}
\put(6,0){\circle*{0.1}}
\put(9,0){\circle*{0.1}}

\put(3,2){\circle*{0.1}}
\put(6,2){\circle*{0.1}}
\put(9,2){\circle*{0.1}}

\put(3,0){\line(0,1){2}}
\put(6,0){\line(-3,2){3}}
\put(6,0){\line(0,1){2}}
\put(9,0){\line(-3,2){3}}
\put(9,0){\line(0,1){2}}
\put(3.1,2.1){$1$}
\put(6.1,2.1){$2$}
\put(9.1,2.1){$3$}
\put(3.1,-0.3){$1'$}
\put(6.1,-0.3){$2'$}
\put(9.1,-0.3){$3'$}
\end{picture}

By Proposition \ref{thm:domaintrees}, $\AG$ is not a domain.

It is easy to check that $G$ has five different basic $1$-covers. Moreover, the number of basic $k$-covers $\aa=(a_1,a_2,a_3,a_{1'},a_{2'},a_{3'})$ such that $\# \{a_1,a_2,a_3\}\leq 2$ is asymptotically smaller than a constant times $k^2$. Notice also that if $\# \{a_1,a_2,a_3\}=3$, then one necessarily has $a_1>a_2>a_3$. From this observations one can deduce that the number of basic $k$-covers of $G$ grows asymptotically as $\displaystyle \frac{1}{6} k^3$.

So $\dim(\AG)=3+1=4$ and $e(\AG)=1$. The height of the ideal $J$ which presents $\AG$ as a quotient of $P$ is 1, since $G$ has five basic 1-covers. Concluding
\[e(\AG)=1<2=1+\operatorname{height}(J),\]
which, together to the trivial fact that $J\subseteq \mm^2$, implies that $\AG$ is not Cohen-Macaulay.
\end{es}

\begin{os}
Actually, in Example \ref{G not CM}, $\AG$ is not even equidimensional; were it equidimensional, $\AG$ would be a complete intersection, which is stronger than the Cohen-Macaulay property.

Example \ref{G not CM} shows also that the two integers $\gdim(G)$ and $\depth(\AG)$, which both are smaller or equal than $\dim \AG$, are in general different.
\end{os}

\vspace{2mm}
%

\subsection{The dimension of the fiber cone associated to a tree}\label{tree}

We have already seen a few graphs for which the graphical dimension of $G$ equals the dimension of the algebra $\AG$. In this subsection we will prove that equality holds also if $G$ is a tree (Theorem \ref{6}). The proof is quite technical, the main idea being to find a certain optimal standard drawing of $G$. We will show that in every optimal standard drawing $\Gamma$ of a connected bipartite graph $G$, if there is an edge to the right of $\{ r(\Gamma), \, r(\Gamma)'\}$, then $G$ must contain a cycle (Proposition \ref{cycle}).

\begin{defn}\label{2}
Fix a standard drawing $\Gamma$ of a bipartite connected graph $G$. Two vertices $i \leq j$ in $A$ are \textit{saw-connected} (with respect to $\Gamma$) if:
\begin{compactitem}[--]
\item $i, j$ are both smaller or equal than $r \left(\Gamma \right)$, and
\item $\Gamma$ contains the backslash edges $\{i_1, i_2'\}$,
$\{i_2, i_3'\} \ldots ,\{i_{t-1}, i_t'\}$, for some integers $i = i_1 <  \ldots < i_t = j$.
 \end{compactitem}
 \end{defn}

We will need the following notation for a connected bipartite graph $G$ and a drawing $\Gamma$.
\begin{compactitem}
\item For every vertex $v \in A$, \  $\mathfrak{j}(v,\Gamma)\in \{1, \ldots, a\}$ will be the position of $v$ in the drawing $\Gamma$;
\item for every $w \in B$, the integer $\mathfrak{j}'(w,\Gamma)\in \{1,\ldots ,b\}$ is defined analogously;
\item for all $i =1, \ldots ,a$ \ \ $\mathfrak{v}(i,\Gamma)\in V(G)$ will be the vertex corresponding to $i$ in the drawing $\Gamma$;
\item for every $j=1, \ldots ,b$ the integer $\vv(j' , \Gamma)$ is defined analogously.
\end{compactitem}

The above notation avoids confusions when passing from a drawing to another.

\begin{lem}\label{saw}
Let $\Gamma$ be a standard drawing of $G$ and let $r=r(\Gamma)$.
\begin{compactenum}
\item Let $\{i,j'\} \in E$ with $i,j>r$: If $\Gamma$ is optimal then
\[m(i,\Gamma) \leq M(j',\Gamma). \]
\item If $i, j$ are not saw-connected with $i<j\leq r$, then there is another standard drawing $\Gamma'$ such that:
    \begin{compactitem}
     \item $\jj(\vv(i,\Gamma), \Gamma')=\jj'(\vv(i',\Gamma), \Gamma')=\jj(\vv(j,\Gamma), \Gamma')+1=\jj'(\vv(j',\Gamma), \Gamma')+1$;
    \item $\Gamma$ and $\Gamma'$ coincide outside the area delimited by the edges $\{i, i'\}$ and $\{j, j'\}$;
    \item $r(\Gamma')=r(\Gamma)$.
    \end{compactitem}
\end{compactenum}
\end{lem}
\begin{proof}
(1). By contradiction assume that $m(i, \Gamma)>M(j',\Gamma)$. Without loss of generality we can also assume $i=j$.\\
Look at the picture below: in that situation $i=j=5$, $m(i,\Gamma)=3$, $M(j',\Gamma)=2$ and $r=3$. We only have to move the edge $\{i,i'\}$ to the $(m(i,\Gamma)+1)$-th position and to move to one step to the right every edge $\{h,h'\}$ for all $h=m(i, \Gamma)+1, \ldots , i-1$. It is easy to see that  $r(\Gamma')=r(\Gamma)+1$, and this contradicts the fact that $\Gamma$ is optimal.

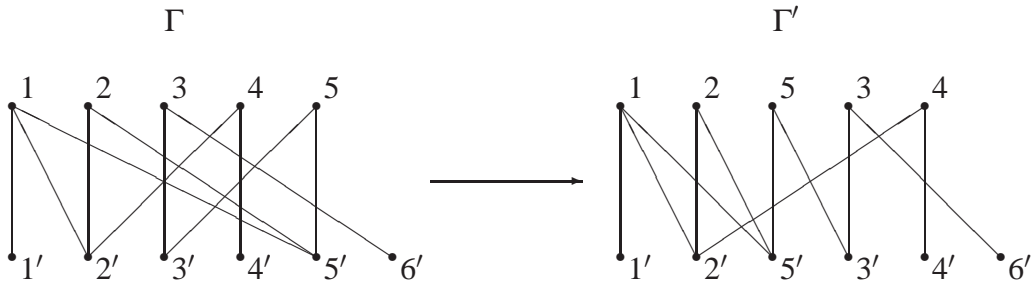
\setlength{\unitlength}{1cm}
\begin{picture}(15,5)(0,-1)
\put(2,3){$\Gamma$}
\put(10,3){$\Gamma'$}

\put(0,0){\circle*{0.1}}
\put(1,0){\circle*{0.1}}
\put(2,0){\circle*{0.1}}
\put(3,0){\circle*{0.1}}
\put(4,0){\circle*{0.1}}
\put(5,0){\circle*{0.1}}

\put(0,2){\circle*{0.1}}
\put(1,2){\circle*{0.1}}
\put(2,2){\circle*{0.1}}
\put(3,2){\circle*{0.1}}
\put(4,2){\circle*{0.1}}

\put(0,0){\line(0,1){2}}
\put(1,0){\line(-1,2){1}}
\put(1,0){\line(0,1){2}}
\put(1,0){\line(1,1){2}}
\put(2,0){\line(1,1){2}}
\put(2,0){\line(0,1){2}}
\put(3,0){\line(0,1){2}}
\put(4,0){\line(-2,1){4}}
\put(4,0){\line(-3,2){3}}
\put(4,0){\line(0,1){2}}
\put(5,0){\line(-3,2){3}}
\put(0.1,2.1){$1$}
\put(1.1,2.1){$2$}
\put(2.1,2.1){$3$}
\put(3.1,2.1){$4$}
\put(4.1,2.1){$5$}
\put(0.1,-0.3){$1'$}
\put(1.1,-0.3){$2'$}
\put(2.1,-0.3){$3'$}
\put(3.1,-0.3){$4'$}
\put(4.1,-0.3){$5'$}
\put(5.1,-0.3){$6'$}

\put(5.5,1){\vector(1,0){2}}

\put(8,0){\circle*{0.1}}
\put(9,0){\circle*{0.1}}
\put(10,0){\circle*{0.1}}
\put(11,0){\circle*{0.1}}
\put(12,0){\circle*{0.1}}
\put(13,0){\circle*{0.1}}

\put(8,2){\circle*{0.1}}
\put(9,2){\circle*{0.1}}
\put(10,2){\circle*{0.1}}
\put(11,2){\circle*{0.1}}
\put(12,2){\circle*{0.1}}

\put(8,0){\line(0,1){2}}
\put(9,0){\line(-1,2){1}}
\put(9,0){\line(0,1){2}}
\put(9,0){\line(3,2){3}}
\put(10,0){\line(-1,2){1}}
\put(10,0){\line(0,1){2}}
\put(10,0){\line(-1,1){2}}
\put(11,0){\line(-1,2){1}}
\put(11,0){\line(0,1){2}}
\put(12,0){\line(0,1){2}}
\put(13,0){\line(-1,1){2}}
\put(8.1,2.1){$1$}
\put(9.1,2.1){$2$}
\put(10.1,2.1){$5$}
\put(11.1,2.1){$3$}
\put(12.1,2.1){$4$}
\put(8.1,-0.3){$1'$}
\put(9.1,-0.3){$2'$}
\put(10.1,-0.3){$5'$}
\put(11.1,-0.3){$3'$}
\put(12.1,-0.3){$4'$}
\put(13.1,-0.3){$6'$}

\end{picture}

\vspace{1mm}

(2). Consider the two sets:
\[C=\{t \in \{i,i+1, \ldots ,j-1,j\}: i,t \mbox{ saw-connected}\}\]
\[N=\{i,i+1, \ldots ,j-1,j\}\setminus C\]
Notice that $i\in C$ and $j \in N$. Write
\[C=\{i_1, \ldots, i_d\}, \mbox{ with }i=i_1 \leq \ldots \leq i_d < j, \hbox{ and}\]
\[N=\{j_1, \ldots, j_v\}, \mbox{ with }i<j_0 \leq \ldots \leq j_v= j,\]
where obviously $v+d=j-i$.

Consider the new drawing $\Gamma'$ obtained by $\Gamma$ using the permutations $\sigma \in S_n$ and $\tt \in S_m$ defined as follows.
\begin{displaymath}
\sigma(k) = \left\{ \begin{array}{cc} k  & \mbox{if } k < i \mbox{  or if }k>j\\
i_l  & \mbox{if } k=i+v+l  \\
j_h  & \mbox{if } k=i+h \mbox{ with } h \leq v \end{array}  \right.
\end{displaymath}
\begin{displaymath}
\tt(k) = \left\{ \begin{array}{cc} k  & \mbox{if } k < i \mbox{  or if }k>j \\
i_l  & \mbox{if } k=i+v+l \\
j_h  & \mbox{if } k=i+h \mbox{ with } h \leq v \end{array}  \right.
\end{displaymath}
It is easy to see that $r(\Gamma)=r(\Gamma')$.
\end{proof}

Look at the picture below for an example of how the proof of point (2) of the lemma above works. In that situation $i=2$, $j=7$, $C=\{2,4,5\}$, $N=\{3,6,7\}$ and $r=7$. We have to move all the edges $\{h,h'\}$ with $h \in C$ to the right of $\{j,j'\}$, without changing their order.

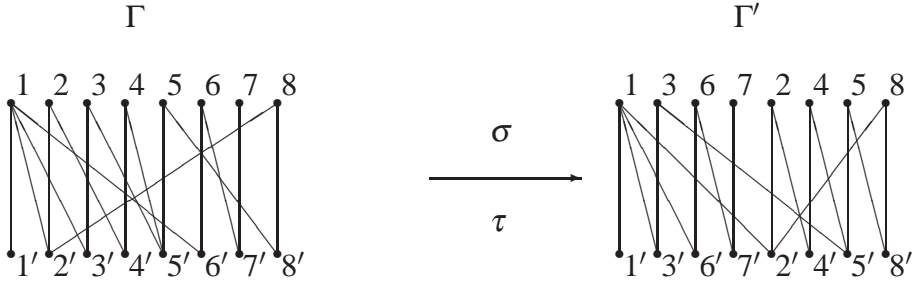
\setlength{\unitlength}{1cm}
\begin{picture}(15,5)(0,-1)
\put(1.5,3){$\Gamma$}
\put(9.5,3){$\Gamma'$}

\put(0,0){\circle*{0.1}}
\put(0.5,0){\circle*{0.1}}
\put(1,0){\circle*{0.1}}
\put(1.5,0){\circle*{0.1}}
\put(2,0){\circle*{0.1}}
\put(2.5,0){\circle*{0.1}}
\put(3,0){\circle*{0.1}}
\put(3.5,0){\circle*{0.1}}

\put(0,2){\circle*{0.1}}
\put(0.5,2){\circle*{0.1}}
\put(1,2){\circle*{0.1}}
\put(1.5,2){\circle*{0.1}}
\put(2,2){\circle*{0.1}}
\put(2.5,2){\circle*{0.1}}
\put(3,2){\circle*{0.1}}
\put(3.5,2){\circle*{0.1}}

\put(0,0){\line(0,1){2}}
\put(0.5,0){\line(3,2){3}}
\put(0.5,0){\line(0,1){2}}
\put(0.5,0){\line(-1,4){0.5}}
\put(1,0){\line(0,1){2}}
\put(1,0){\line(-1,2){1}}
\put(1.5,0){\line(-1,2){1}}
\put(1.5,0){\line(0,1){2}}
\put(2,0){\line(-1,4){0.5}}
\put(2,0){\line(-1,2){1}}
\put(2,0){\line(0,1){2}}
\put(2.5,0){\line(0,1){2}}
\put(2.5,0){\line(-5,4){2.5}}
\put(3,0){\line(-1,4){0.5}}
\put(3,0){\line(0,1){2}}
\put(3.5,0){\line(0,1){2}}
\put(3.5,0){\line(-3,4){1.5}}
\put(0.05,2.1){$1$}
\put(0.55,2.1){$2$}
\put(1.05,2.1){$3$}
\put(1.55,2.1){$4$}
\put(2.05,2.1){$5$}
\put(3.05,2.1){$7$}
\put(2.55,2.1){$6$}
\put(3.55,2.1){$8$}
\put(0.05,-0.3){$1'$}
\put(0.55,-0.3){$2'$}
\put(1.05,-0.3){$3'$}
\put(1.55,-0.3){$4'$}
\put(2.05,-0.3){$5'$}
\put(3.05,-0.3){$7'$}
\put(2.55,-0.3){$6'$}
\put(3.55,-0.3){$8'$}

\put(5.5,1){\vector(1,0){2}}
\put(6.3,1.5){$\sigma$}
\put(6.3,0.3){$\tt$}

\put(8,0){\circle*{0.1}}
\put(8.5,0){\circle*{0.1}}
\put(9,0){\circle*{0.1}}
\put(9.5,0){\circle*{0.1}}
\put(10,0){\circle*{0.1}}
\put(10.5,0){\circle*{0.1}}
\put(11,0){\circle*{0.1}}
\put(11.5,0){\circle*{0.1}}

\put(8,2){\circle*{0.1}}
\put(8.5,2){\circle*{0.1}}
\put(9,2){\circle*{0.1}}
\put(9.5,2){\circle*{0.1}}
\put(10,2){\circle*{0.1}}
\put(10.5,2){\circle*{0.1}}
\put(11,2){\circle*{0.1}}
\put(11.5,2){\circle*{0.1}}

\put(8,0){\line(0,1){2}}
\put(10,0){\line(3,4){1.5}}
\put(8.5,0){\line(0,1){2}}
\put(8.5,0){\line(-1,4){0.5}}
\put(9,0){\line(0,1){2}}
\put(9,0){\line(-1,2){1}}
\put(9.5,0){\line(-1,4){0.5}}
\put(9.5,0){\line(0,1){2}}
\put(11,0){\line(-1,4){0.5}}
\put(10,0){\line(-1,1){2}}
\put(10,0){\line(0,1){2}}
\put(10.5,0){\line(0,1){2}}
\put(10.5,0){\line(-1,4){0.5}}
\put(11,0){\line(-5,4){2.5}}
\put(11,0){\line(0,1){2}}
\put(11.5,0){\line(0,1){2}}
\put(11.5,0){\line(-1,4){0.5}}
\put(8.05,2.1){$1$}
\put(8.55,2.1){$3$}
\put(9.05,2.1){$6$}
\put(9.55,2.1){$7$}
\put(10.05,2.1){$2$}
\put(10.55,2.1){$4$}
\put(11.05,2.1){$5$}
\put(11.55,2.1){$8$}
\put(8.05,-0.3){$1'$}
\put(8.55,-0.3){$3'$}
\put(9.05,-0.3){$6'$}
\put(9.55,-0.3){$7'$}
\put(10.05,-0.3){$2'$}
\put(10.55,-0.3){$4'$}
\put(11.05,-0.3){$5'$}
\put(11.55,-0.3){$8'$}

\end{picture}

\begin{prop}\label{cycle}
Let $G$ be a connected bipartite graph and $\Gamma$ an optimal standard drawing. If $\{i,j' \} \in E$ for $i,j>r$ then $\{i,j'\}$ is part of a cycle.
\end{prop}
\begin{proof}
By point (1) of Lemma \ref{saw} $m(i,\Gamma) \leq M(j',\Gamma)$. We prove the statement by induction on $s=M(j',\Gamma)-m(i,\Gamma)$.\\
When $s=0$, set $z:= m(i,\Gamma)=M(j',\Gamma)$. Since $z \leq r$, there is a vertical edge going from $z$ in $A$ to $z'$ in $B$; so $G$ contains the cycle $\{z, z', i, j'\}$.\\
Let $s$ be bigger than $0$. If $\{i,j'\}$ is not a part of a cycle then $m(i,\Gamma)$ and $M(j',\Gamma)$ are not saw-connected. By point (2) of Lemma \ref{saw} there is another standard drawing $\Gamma'$ such that $\jj(\vv(M(j',\Gamma),\Gamma')<\jj(\vv(M(i,\Gamma),\Gamma')$, $\jj(\vv(i,\Gamma), \Gamma')=i$ and $\jj(\vv(j',\Gamma),\Gamma')=j$ . It is clear that $M(j',\Gamma')<M(j',\Gamma)$ (otherwise there would be a cycle), and that $m(i,\Gamma')\geq m(i,\Gamma)$, so we conclude by induction.
\end{proof}

\begin{prop}\label{4}
Suppose that $G$ is a tree. Then there is an optimal standard drawing $\Gamma$ of $G$ such that to the right of the edge $\{r(\Gamma),r(\Gamma)'\}$ there are only leaves.
\end{prop}
\begin{proof}
During the proof, for every optimal drawing $\Gamma$ and for every $j>r(\Gamma)$ such that $j'$ is not a leaf, we will denote by
\[N(j',\Gamma):=\max\{i\neq M(j',\Gamma):\{i,j'\} \mbox{ is an edge}\}.\]
Notice that by Proposition \ref{cycle} $N(j',\Gamma)<M(j',\Gamma)\leq r(\Gamma)$.

Consider an optimal  drawing $\Gamma_0$ of $G$ such that the number of non-leaves to the right of $r:=r(\Gamma_0)$ is minimal among the other optimal standard drawings. By contradiction assume there are some non-leaves to the right of $\{ r, r'\}$.

First consider the case in which there is a non-leaf to the right of $\{ r, r'\}$ in $B$. Pick $j>r$ such that $j'$ is not a leaf.
\begin{compactenum}
\item Consider the new, clearly optimal, standard drawing $\Gamma_1$ obtained by switching $M(j',\Gamma_0)'$ and $j'$. If $\vv(j',\Gamma_1)=\vv(M(j',\Gamma_0)', \Gamma_0)$ is a leaf we conclude, because in the optimal drawing $\Gamma_1$ there are less leaves to the right of $r$ than in $\Gamma_0$, and this is a contradiction.
\item Otherwise we can consider $M(j', \Gamma_1)$ and $N(j', \Gamma_1)$. Since $G$ is a tree, $N(j',\Gamma_1)$ and $M(j',\Gamma_1)$ are not saw connected, otherwise  there would be a cycle in $G$; so we can find a new drawing $\Gamma_2$ of $G$ as in Lemma 3.11, point (2). Notice that $\vv(j',\Gamma_2)=\vv(j',\Gamma_1) \neq \vv(j',\Gamma_0)$.
\item Consider the new optimal standard drawing $\Gamma_3$ obtained switching $M(j',\Gamma_2)'$ and $j'$. If $\vv(j',\Gamma_3)=\vv(M(j',\Gamma_2)',\Gamma_2)$ is a leaf we conclude as in point (1).
\item Otherwise, arguing as in point (2), $N(j',\Gamma_3)$ and $M(j',\Gamma_3)$ are not saw connected; so we can find a new drawing $\Gamma_4$ of $G$ as in Lemma 3.11, point (2).
Notice that $\vv(j',\Gamma_4)=\vv(j',\Gamma_3) \neq \vv(j',\Gamma_2)$.\\
Moreover, $\vv(j',\Gamma_4)\neq \vv(j',\Gamma_0)$. In fact one has
\[\vv(j',\Gamma_0)\neq \vv(M(j',\Gamma_0)',\Gamma_0)\]
and
\[ \vv(M(j' ,\Gamma_2), \Gamma_2)=\vv(N(j',\Gamma_1),\Gamma_1)\neq \vv(M(j',\Gamma_1), \Gamma_1)= \vv(M(j' , \Gamma_0), \Gamma_0).\]
Since $\vv(M(j',\Gamma_2)',\Gamma_2)=\vv(j',\Gamma_4)$ it follows  that $\{\vv(M(j' ,\Gamma_2), \Gamma_2), \vv(j' ,\Gamma_4) \}$ is an edge. So if $\vv(j', \Gamma_4)=\vv(j', \Gamma_0)$. But then
\[ \vv(j',\Gamma_0), \ \ \vv(M(j',\Gamma_0),\Gamma_0), \ \ \vv(M(j', \Gamma_0)',\Gamma_0), \ \ \vv(M(j',\Gamma_2),\Gamma_2), \ \ \vv(j', \Gamma_4)=\vv(j', \Gamma_0)\]
would be a cycle, which is a contradiction.
\item Consider the new optimal standard drawing $\Gamma_5$ obtained switching $M(j',\Gamma_4)'$ and $j'$. If $\vv(j',\Gamma_5)$ is a leaf we conclude as in point (1).
\end{compactenum}
and so on. Eventually we obtain a sequence of optimal drawings $(\Gamma_l)_{l\geq 0}$ such that:
\begin{compactitem}[---]
\item $\{\vv(j',\Gamma_{2l}),\vv(M(j',\Gamma_{2l}),\Gamma_{2l})\}$, \ $\{\vv(M(j',\Gamma_{2l}),\Gamma_{2l}), \vv(j',\Gamma_{2l+2})\}$ are edges for every $l\geq 0$.
\item $\vv(j',\Gamma_{2l})\neq \vv(j',\Gamma_{2l+2})\neq \vv(j', \Gamma_{2l+4})\neq \vv(j',\Gamma_{2l})$ for every $l \geq 0$.
\item $\vv(M(j',\Gamma_{2l}),\Gamma_2)$ is not a leaf for any $l\geq 0$.
\end{compactitem}
So starting from $j'$ we can construct stepwise the following path in $G$:
\[\vv(j', \Gamma_0), \ \vv(M(j',\Gamma_0),\Gamma_0), \ \vv(j',\Gamma_2), \ \vv(M(j',\Gamma_2),\Gamma_2), \ \vv(j',\Gamma_4), \ \vv(M(j',\Gamma_4),\Gamma_4), \ \vv(j',\Gamma_6), \ldots \]
As $G$ is a tree and as $\vv(M(j',\Gamma_{2l}),\Gamma_{2l-1})$ is never a leaf, then, at some point, there will be a $l_0$ such that $\vv(j',\Gamma_{2l_0})=\vv(j',\Gamma_{2l_0-1})$ is a leaf. This contradicts the fact that for each $l\geq 0$ we have $\vv(i,\Gamma_l)=\vv(i,\Gamma_0)$ for every $i>r$, and $\vv(k',\Gamma_l)=\vv(k',\Gamma_0)$ for every $k>r$ such that $k \neq j$. So $\Gamma_{2l_0-1}$ has less leaves to right of $\{ r, r'\}$ than $\Gamma_0$.

If, relatively to the drawing $\Gamma_0$, for every $j>r$ \ $j'$ is a leaf, then there exists $i>r$ such that $i$ is not a leaf. The argument is the same of the first part with the only difference that we have to consider $m(i,\Gamma_0)$ and
\[ n(i,\Gamma_0):=\min \{j\neq m(i,\Gamma_0):\{i,j'\} \mbox{ is an edge}\} \]
Notice that $m(i,\Gamma_0)<n(i,\Gamma_0)\leq r$ by Proposition \ref{cycle}.
\end{proof}

\begin{es}\label{alg}
We will show with an example how the proof of the above proposition works.

Below $j=6, \ \vv(j',\Gamma_0)=v_6, \ M(j',\Gamma_0)=4,  \ \vv(M(j', \Gamma_0),\Gamma_0)=u_4,$ $\vv(j', \Gamma_2)=v_4$.

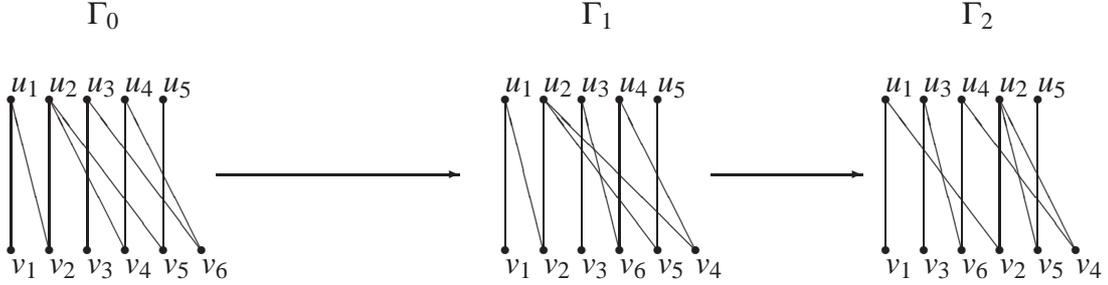
\setlength{\unitlength}{1cm}
\begin{picture}(15,5)(0,-1)
\put(1,3){$\Gamma_0$}
\put(7.5,3){$\Gamma_1$}
\put(12.5,3){$\Gamma_2$}

\put(0,0){\circle*{0.1}}
\put(0.5,0){\circle*{0.1}}
\put(1,0){\circle*{0.1}}
\put(1.5,0){\circle*{0.1}}
\put(2,0){\circle*{0.1}}
\put(2.5,0){\circle*{0.1}}

\put(0,2){\circle*{0.1}}
\put(0.5,2){\circle*{0.1}}
\put(1,2){\circle*{0.1}}
\put(1.5,2){\circle*{0.1}}
\put(2,2){\circle*{0.1}}

\put(0,0){\line(0,1){2}}
\put(0.5,0){\line(0,1){2}}
\put(0.5,0){\line(-1,4){0.5}}
\put(1,0){\line(0,1){2}}
\put(1.5,0){\line(-1,2){1}}
\put(1.5,0){\line(0,1){2}}
\put(2,0){\line(0,1){2}}
\put(2,0){\line(-3,4){1.5}}
\put(2.5,0){\line(-1,2){1}}
\put(2.5,0){\line(-3,4){1.5}}
\put(0,2.1){$u_1$}
\put(0.5,2.1){$u_2$}
\put(1,2.1){$u_3$}
\put(1.5,2.1){$u_4$}
\put(2,2.1){$u_5$}
\put(0,-0.3){$v_1$}
\put(0.5,-0.3){$v_2$}
\put(1,-0.3){$v_3$}
\put(1.5,-0.3){$v_4$}
\put(2,-0.3){$v_5$}
\put(2.5,-0.3){$v_6$}

\put(2.7,1){\vector(1,0){3.2}}

\put(6.5,0){\circle*{0.1}}
\put(7,0){\circle*{0.1}}
\put(7.5,0){\circle*{0.1}}
\put(8,0){\circle*{0.1}}
\put(8.5,0){\circle*{0.1}}
\put(9,0){\circle*{0.1}}

\put(6.5,2){\circle*{0.1}}
\put(7,2){\circle*{0.1}}
\put(7.5,2){\circle*{0.1}}
\put(8,2){\circle*{0.1}}
\put(8.5,2){\circle*{0.1}}

\put(6.5,0){\line(0,1){2}}
\put(7,0){\line(0,1){2}}
\put(7,0){\line(-1,4){0.5}}
\put(7.5,0){\line(0,1){2}}
\put(8,0){\line(-1,4){0.5}}
\put(8,0){\line(0,1){2}}
\put(8.5,0){\line(0,1){2}}
\put(8.5,0){\line(-3,4){1.5}}
\put(9,0){\line(-1,2){1}}
\put(9,0){\line(-1,1){2}}
\put(6.5,2.1){$u_1$}
\put(7,2.1){$u_2$}
\put(7.5,2.1){$u_3$}
\put(8,2.1){$u_4$}
\put(8.5,2.1){$u_5$}
\put(6.5,-0.3){$v_1$}
\put(7,-0.3){$v_2$}
\put(7.5,-0.3){$v_3$}
\put(8,-0.3){$v_6$}
\put(8.5,-0.3){$v_5$}
\put(9,-0.3){$v_4$}

\put(9.2,1){\vector(1,0){2}}

\put(11.5,0){\circle*{0.1}}
\put(12,0){\circle*{0.1}}
\put(12.5,0){\circle*{0.1}}
\put(13,0){\circle*{0.1}}
\put(13.5,0){\circle*{0.1}}
\put(14,0){\circle*{0.1}}

\put(11.5,2){\circle*{0.1}}
\put(12,2){\circle*{0.1}}
\put(12.5,2){\circle*{0.1}}
\put(13,2){\circle*{0.1}}
\put(13.5,2){\circle*{0.1}}

\put(11.5,0){\line(0,1){2}}
\put(12,0){\line(0,1){2}}
\put(12.5,0){\line(-1,4){0.5}}
\put(12.5,0){\line(0,1){2}}
\put(13,0){\line(-3,4){1.5}}
\put(13,0){\line(0,1){2}}
\put(13.5,0){\line(0,1){2}}
\put(13.5,0){\line(-1,4){0.5}}
\put(14,0){\line(-1,2){1}}
\put(14,0){\line(-3,4){1.5}}
\put(11.5,2.1){$u_1$}
\put(12,2.1){$u_3$}
\put(12.5,2.1){$u_4$}
\put(13,2.1){$u_2$}
\put(13.5,2.1){$u_5$}
\put(11.5,-0.3){$v_1$}
\put(12,-0.3){$v_3$}
\put(12.5,-0.3){$v_6$}
\put(13,-0.3){$v_2$}
\put(13.5,-0.3){$v_5$}
\put(14,-0.3){$v_4$}

\end{picture}

Below $ M(j',\Gamma_2)=4, \ \vv(M(j', \Gamma_2),\Gamma_2)=u_2, \ \vv(j', \Gamma_4)=v_2$.

\setlength{\unitlength}{1cm}
\begin{picture}(15,5)(0,-1)
\put(1,3){$\Gamma_2$}
\put(7.5,3){$\Gamma_3$}
\put(12.5,3){$\Gamma_4$}

\put(0,0){\circle*{0.1}}
\put(0.5,0){\circle*{0.1}}
\put(1,0){\circle*{0.1}}
\put(1.5,0){\circle*{0.1}}
\put(2,0){\circle*{0.1}}
\put(2.5,0){\circle*{0.1}}

\put(0,2){\circle*{0.1}}
\put(0.5,2){\circle*{0.1}}
\put(1,2){\circle*{0.1}}
\put(1.5,2){\circle*{0.1}}
\put(2,2){\circle*{0.1}}

\put(11.5,0){\circle*{0.1}}
\put(12,0){\circle*{0.1}}
\put(12.5,0){\circle*{0.1}}
\put(13,0){\circle*{0.1}}
\put(13.5,0){\circle*{0.1}}
\put(14,0){\circle*{0.1}}

\put(11.5,2){\circle*{0.1}}
\put(12,2){\circle*{0.1}}
\put(12.5,2){\circle*{0.1}}
\put(13,2){\circle*{0.1}}
\put(13.5,2){\circle*{0.1}}

\put(6.5,0){\circle*{0.1}}
\put(7,0){\circle*{0.1}}
\put(7.5,0){\circle*{0.1}}
\put(8,0){\circle*{0.1}}
\put(8.5,0){\circle*{0.1}}
\put(9,0){\circle*{0.1}}

\put(6.5,2){\circle*{0.1}}
\put(7,2){\circle*{0.1}}
\put(7.5,2){\circle*{0.1}}
\put(8,2){\circle*{0.1}}
\put(8.5,2){\circle*{0.1}}

\put(0,0){\line(0,1){2}}
\put(0.5,0){\line(0,1){2}}
\put(1,0){\line(-1,4){0.5}}
\put(1,0){\line(0,1){2}}
\put(1.5,0){\line(-3,4){1.5}}
\put(1.5,0){\line(0,1){2}}
\put(2,0){\line(0,1){2}}
\put(2,0){\line(-1,4){0.5}}
\put(2.5,0){\line(-1,2){1}}
\put(2.5,0){\line(-3,4){1.5}}
\put(0,2.1){$u_1$}
\put(0.5,2.1){$u_3$}
\put(1,2.1){$u_4$}
\put(1.5,2.1){$u_2$}
\put(2,2.1){$u_5$}
\put(0,-0.3){$v_1$}
\put(0.5,-0.3){$v_3$}
\put(1,-0.3){$v_6$}
\put(1.5,-0.3){$v_2$}
\put(2,-0.3){$v_5$}
\put(2.5,-0.3){$v_4$}

\put(2.7,1){\vector(1,0){3.2}}

\put(6.5,0){\line(0,1){2}}
\put(7,0){\line(0,1){2}}
\put(7.5,0){\line(-1,4){0.5}}
\put(7.5,0){\line(0,1){2}}
\put(8,0){\line(-1,4){0.5}}
\put(8,0){\line(0,1){2}}
\put(8.5,0){\line(0,1){2}}
\put(8.5,0){\line(-1,4){0.5}}
\put(9,0){\line(-1,2){1}}
\put(9,0){\line(-5,4){2.5}}
\put(6.5,2.1){$u_1$}
\put(7,2.1){$u_3$}
\put(7.5,2.1){$u_4$}
\put(8,2.1){$u_2$}
\put(8.5,2.1){$u_5$}
\put(6.5,-0.3){$v_1$}
\put(7,-0.3){$v_3$}
\put(7.5,-0.3){$v_6$}
\put(8,-0.3){$v_4$}
\put(8.5,-0.3){$v_5$}
\put(9,-0.3){$v_2$}

\put(9.2,1){\vector(1,0){2}}

\put(11.5,0){\line(0,1){2}}
\put(12,0){\line(0,1){2}}
\put(12,0){\line(-1,4){0.5}}
\put(12.5,0){\line(0,1){2}}
\put(12.5,0){\line(-1,4){0.5}}
\put(13,0){\line(0,1){2}}
\put(13.5,0){\line(0,1){2}}
\put(13.5,0){\line(-1,2){1}}
\put(14,0){\line(-1,2){1}}
\put(14,0){\line(-3,4){1.5}}
\put(11.5,2.1){$u_3$}
\put(12,2.1){$u_4$}
\put(12.5,2.1){$u_2$}
\put(13,2.1){$u_1$}
\put(13.5,2.1){$u_5$}
\put(11.5,-0.3){$v_3$}
\put(12,-0.3){$v_6$}
\put(12.5,-0.3){$v_4$}
\put(13,-0.3){$v_1$}
\put(13.5,-0.3){$v_5$}
\put(14,-0.3){$v_2$}

\end{picture}

Below $M(j',\Gamma_4)=4, \ \vv(M(j',\Gamma_4),\Gamma_4)=u_1, \  \vv(j', \Gamma_5)=v_1$.

\setlength{\unitlength}{1cm}
\begin{picture}(15,5)(0,-1)
\put(1,3){$\Gamma_4$}
\put(7.5,3){$\Gamma_5$}

\put(0,0){\circle*{0.1}}
\put(0.5,0){\circle*{0.1}}
\put(1,0){\circle*{0.1}}
\put(1.5,0){\circle*{0.1}}
\put(2,0){\circle*{0.1}}
\put(2.5,0){\circle*{0.1}}

\put(0,2){\circle*{0.1}}
\put(0.5,2){\circle*{0.1}}
\put(1,2){\circle*{0.1}}
\put(1.5,2){\circle*{0.1}}
\put(2,2){\circle*{0.1}}

\put(6.5,0){\circle*{0.1}}
\put(7,0){\circle*{0.1}}
\put(7.5,0){\circle*{0.1}}
\put(8,0){\circle*{0.1}}
\put(8.5,0){\circle*{0.1}}
\put(9,0){\circle*{0.1}}

\put(6.5,2){\circle*{0.1}}
\put(7,2){\circle*{0.1}}
\put(7.5,2){\circle*{0.1}}
\put(8,2){\circle*{0.1}}
\put(8.5,2){\circle*{0.1}}

\put(0,0){\line(0,1){2}}
\put(0.5,0){\line(0,1){2}}
\put(0.5,0){\line(-1,4){0.5}}
\put(1,0){\line(0,1){2}}
\put(1,0){\line(-1,4){0.5}}
\put(1.5,0){\line(0,1){2}}
\put(2,0){\line(0,1){2}}
\put(2,0){\line(-1,2){1}}
\put(2.5,0){\line(-1,2){1}}
\put(2.5,0){\line(-3,4){1.5}}
\put(0,2.1){$u_3$}
\put(0.5,2.1){$u_4$}
\put(1,2.1){$u_2$}
\put(1.5,2.1){$u_1$}
\put(2,2.1){$u_5$}
\put(0,-0.3){$v_3$}
\put(0.5,-0.3){$v_6$}
\put(1,-0.3){$v_4$}
\put(1.5,-0.3){$v_1$}
\put(2,-0.3){$v_5$}
\put(2.5,-0.3){$v_2$}

\put(2.7,1){\vector(1,0){3.2}}

\put(6.5,0){\line(0,1){2}}
\put(7,0){\line(0,1){2}}
\put(7,0){\line(-1,4){0.5}}
\put(7.5,0){\line(0,1){2}}
\put(7.5,0){\line(-1,4){0.5}}
\put(8,0){\line(0,1){2}}
\put(8.5,0){\line(0,1){2}}
\put(8,0){\line(-1,4){0.5}}
\put(8.5,0){\line(-1,2){1}}
\put(9,0){\line(-1,2){1}}
\put(6.5,2.1){$u_3$}
\put(7,2.1){$u_4$}
\put(7.5,2.1){$u_2$}
\put(8,2.1){$u_1$}
\put(8.5,2.1){$u_5$}
\put(6.5,-0.3){$v_3$}
\put(7,-0.3){$v_6$}
\put(7.5,-0.3){$v_4$}
\put(8,-0.3){$v_2$}
\put(8.5,-0.3){$v_5$}
\put(9,-0.3){$v_1$}

\end{picture}

So we have constructed the path
\[v_6, \ u_4, \ v_4, \ u_2, \ v_2, \ u_1, \ v_1\]

\end{es}

\begin{lem}\label{5}
Let $G$ be a tree and $\Gamma$ an optimal standard drawing of $G$ like that in Proposition \ref{4}. Assume that for some $j< \gdim(G)$, $\{i,j'\}$ is an edge for an $i\geq \gdim(G)$.
If $\{k,h'\}$ is an edge for some $k< \gdim(G)$ and a $h\geq \gdim(G)$, then $j$ and $k$ are not saw-connected.
\end{lem}
\begin{proof}
By contradiction, let $j\leq k < \gdim(G)$ be such that $j$ and $k$ are saw-connected and $\{k,h'\}$ is an edge for some $h>r$. Consider a sequence $j=i_1 < i_2  < \ldots < i_s=k$ such that $\{i_q,i_{q+1}'\}$ is an edge for every $q=1, \ldots, s-1$. Since $G$ is acyclic, this sequence is unique. Applying several times point (2) of Lemma \ref{saw}, we can find another optimal standard drawing $\Gamma'$ for which $i_{q+1}=i_q+1$ for all $q \in \{1, \ldots ,s-1\}$.

We can assume, without loss of generality, that $i_{q+1}=i_q+1$ for each $q$ in $\{1, \ldots ,s-1\}$; we can also assume $i=h$.
Then we can pass to a new drawing $\Gamma'$ using the following permutations:
\begin{displaymath}
\sigma(l) = \left\{ \begin{array}{ll} l  & \mbox{if } l< j \mbox{  or if }l>i \\
i_s-q+1 &  \mbox{if } l=i_q \\
i  & \mbox{if } l=k+1 \\
l-1 & \mbox{if } k+1<l\leq i
\end{array}  \right.
\end{displaymath}
\begin{displaymath}
\tt(l) = \left\{ \begin{array}{ll} l  & \mbox{if } l< j \mbox{  or if }l>i \\
i_s-q+2 &  \mbox{if } l=i_q, q\geq 2 \\
i  & \mbox{if } l=j \\
j & \mbox{if } l=k+1 \\
l-1 & \mbox{if } k+1<l\leq i
\end{array}  \right.
\end{displaymath}
Now it is not difficult to verify that $r(\Gamma')=r(\Gamma)+1$, a contradiction.
\end{proof}

We explain with an example the idea of the proof above. In that situation $i=6, \ j=2, \ k=4, \ s=3, \ i_1=2, \ i_2=3, \ i_3=4$.  We do the following moves:
\begin{compactitem}[-----]
\item $6 \rightarrow 5$;
\item $5 \rightarrow 6$;
\item $2' \rightarrow 5'$;
\item ''walk backwards'' on the path that connects $\{2,4\}$;
\item $6' \rightarrow 2'$.
\end{compactitem}
\setlength{\unitlength}{1cm}
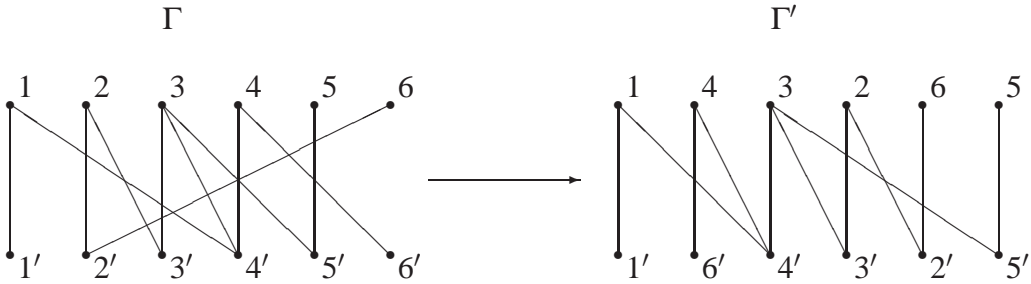
\begin{picture}(15,5)(0,-1)
\put(2,3){$\Gamma$}
\put(10,3){$\Gamma'$}

\put(0,0){\circle*{0.1}}
\put(1,0){\circle*{0.1}}
\put(2,0){\circle*{0.1}}
\put(3,0){\circle*{0.1}}
\put(4,0){\circle*{0.1}}
\put(5,0){\circle*{0.1}}

\put(0,2){\circle*{0.1}}
\put(1,2){\circle*{0.1}}
\put(2,2){\circle*{0.1}}
\put(3,2){\circle*{0.1}}
\put(4,2){\circle*{0.1}}
\put(5,2){\circle*{0.1}}

\put(8,0){\circle*{0.1}}
\put(9,0){\circle*{0.1}}
\put(10,0){\circle*{0.1}}
\put(11,0){\circle*{0.1}}
\put(12,0){\circle*{0.1}}
\put(13,0){\circle*{0.1}}

\put(8,2){\circle*{0.1}}
\put(9,2){\circle*{0.1}}
\put(10,2){\circle*{0.1}}
\put(11,2){\circle*{0.1}}
\put(12,2){\circle*{0.1}}
\put(13,2){\circle*{0.1}}

\put(0,0){\line(0,1){2}}
\put(2,0){\line(-1,2){1}}
\put(1,0){\line(0,1){2}}
\put(1,0){\line(2,1){4}}
\put(2,0){\line(0,1){2}}
\put(3,0){\line(0,1){2}}
\put(3,0){\line(-3,2){3}}
\put(3,0){\line(-1,2){1}}
\put(4,0){\line(-1,1){2}}
\put(4,0){\line(0,1){2}}
\put(5,0){\line(-1,1){2}}
\put(0.1,2.1){$1$}
\put(1.1,2.1){$2$}
\put(2.1,2.1){$3$}
\put(3.1,2.1){$4$}
\put(4.1,2.1){$5$}
\put(5.1,2.1){$6$}
\put(0.1,-0.3){$1'$}
\put(1.1,-0.3){$2'$}
\put(2.1,-0.3){$3'$}
\put(3.1,-0.3){$4'$}
\put(4.1,-0.3){$5'$}
\put(5.1,-0.3){$6'$}

\put(5.5,1){\vector(1,0){2}}

\put(8,0){\line(0,1){2}}
\put(9,0){\line(0,1){2}}
\put(10,0){\line(-1,1){2}}
\put(10,0){\line(-1,2){1}}
\put(10,0){\line(0,1){2}}
\put(11,0){\line(-1,2){1}}
\put(11,0){\line(0,1){2}}
\put(12,0){\line(0,1){2}}
\put(12,0){\line(-1,2){1}}
\put(13,0){\line(-3,2){3}}
\put(13,0){\line(0,1){2}}
\put(8.1,2.1){$1$}
\put(9.1,2.1){$4$}
\put(10.1,2.1){$3$}
\put(11.1,2.1){$2$}
\put(12.1,2.1){$6$}
\put(13.1,2.1){$5$}
\put(8.1,-0.3){$1'$}
\put(9.1,-0.3){$6'$}
\put(10.1,-0.3){$4'$}
\put(11.1,-0.3){$3'$}
\put(12.1,-0.3){$2'$}
\put(13.1,-0.3){$5'$}

\end{picture}

We are now ready to measure the graphical dimension of every tree.

\begin{thm}\label{6}
Suppose that $G$ is a tree. Then
\[\dim(\AG)=\gdim(G).\]
Moreover, if $r=\gdim(G)-1$,
\[e(\AG) \leq (a-r)^{r} \cdot a \cdot (a-1) \cdots (a- r+1)\]
\end{thm}
\begin{proof}
Consider an optimal standard drawing $\Gamma$ of $G$ as in Lemma \ref{4} and set $r=r(\Gamma)=\gdim(G)-1$. We claim that
\begin{equation}
\left \{
    \begin{array}{c}
    \hbox{basic $k$-covers} \\
    \hbox{of the graph } G
    \end{array}
\right \}
\subseteq
\left \{
    \begin{array}{c}
    \left(v_1, v_2, \ldots, v_a, w_{1'}, w_{2'}, \ldots, w_{b'} \right) \\
    \hbox{ in } \{0, 1, \ldots, k\}^n , \hbox{ with } \# \left\{  v_1, \ldots, v_a \right \} \leq r \\
    \hbox{ and }  w_{j'}=k - \min \left\{ v_i ~|~ \{i,j'\} \hbox{ is an edge of } G \right\}.
    \end{array}
\!\!\!\right \}
\end{equation}
Denote by $J_k$ the set on the right hand side.

Suppose that $\aa$ is a basic $k$-cover with distinct $a_i$'s for $i\leq r$. We have to show that for all $j >r$ in $A$, $a_j$ is equal to $a_i$ for some $i \leq r$. By contradiction, suppose there is a $j \in \{r+1, \ldots, a\}$ such that $a_j \neq a_i$ for every $i=1, \ldots,r$. Since $j>r$, the vertex labeled by $j$ is a leaf (Proposition \ref{4}), so there is a unique $i$ (necessarily smaller or equal than $r$) such that $\{j,i'\}$ is an edge of $G$. Clearly $a_{i'}$ must be equal to $k-a_j$. Since $a_i\neq a_j$ there exists $i_1> i$ such that $a_{i_1'}=k-a_i$. Notice that Lemma \ref{5} implies that $i_1 \leq r$. Then, since $a_{i_1}\neq a_i$, there must be an $i_2>i_1$ such that $a_{i_2'}=k-a_{i_1}$. Again Lemma \ref{5} implies that $i_2\leq r$. Continuing this way we find an infinite sequence of integers
\[1\leq i_1<i_2< \ldots \leq r,\]
a contradiction.\\
If $a_1, \ldots, a_r$ are not all distinct, either there is a $j \in \{r+1, \ldots, a\}$ such that $a_j \neq a_i$ for every $i=1, \ldots,r$, or not. If not, the cardinality of $\{a_1, \ldots, a_a\}$ is strictly smaller than $r$, and we are done. Otherwise, take such a $j > r$ with the property that $a_j \neq a_i$ for every $i=1, \ldots,r$; this $j$ is a leaf in $G$, by Proposition \ref{4}; let $k$ be the unique integer smaller than $r$ such that $\{j,k'\}$ is an edge of $G$. There is a $h$ such that $k\neq h\leq r$ and $a_k=a_h$, otherwise one could find a contradiction as before. Moreover if $\{j,k'\}$ and $\{l,k'\}$ are edges for some $j,l>r$ then $a_j=a_l$, because they are leaves. So
\[\# \{a_1, \ldots, a_a\} \leq r.\]

Finally, $a_{j'}=k-\min\{a_i:\{i,j'\} \hbox{ is an edge of } G\}$ for every $j'\in B$, otherwise $\aa$ would not be a basic $k$-cover. So
\[\{\mbox{basic }k \mbox{-covers of }G \} \subseteq J_k\]
Now it is easy to show that
\[\# J_k = (a-r)^{r} \cdot \binom{a}{r} \cdot k^r + O(k^{r-1}),\]
and this ends the proof.
\end{proof}

We conclude this subsection observing that in the tree case we have sufficient conditions for a standard drawing to be optimal.

\begin{prop}
Let $G$ be a tree and $\Gamma$ a standard drawing of $G$ such that there are only leaves to the right of the edge $\{r(\Gamma),r(\Gamma)'\}$. Then the following are equivalent:
\begin{compactenum}
\item $\Gamma$ is optimal;
\item if $\{i,j'\}$ is an edge for some $i>r(\Gamma)$ and $\{k,h'\}$ is an edge for some $h>r(\Gamma)$, then $j$ and $k$ are not saw-connected.
\end{compactenum}
\end{prop}
\begin{proof}
That (1) $\implies$ (2) follows from Lemma \ref{5}; the converse follows from the proof of Theorem \ref{6} (in fact in that proof we used only property (2) of the chosen drawing).
\end{proof}

\section{Further results}

In this last section we describe some further results concerning the particular case of unmixed graphs, the general case of hypergraphs, and some applications to the arithmetical rank of certain monomial ideals.

\subsection{Unmixed graphs}

In this subsection we discuss some other applications of Theorem \ref{thm:domain} with a particular attention to the unmixed case.
Bertone and Micale, during PRAGMATIC 2008, found an example of an unmixed graph $G$ such that $\AG$ is not Gorenstein. Let $G$ be a Cohen-Macaulay graph; when $\AG$ is Gorenstein? Following a hint of Herzog, we were able to answer the question by proving a more general statement: namely, we will show that $\AG$ is an Hibi ring for any unmixed bipartite graph $G$ (see Theorem \ref{hibiring}).

For our discussion, notice that we can present $\AG$ as a quotient of the polynomial ring $P=\KK[X_{\aa}:\aa  \mbox{ is a basic 1-cover of }G]$, namely
\[P/J \cong \AG \]
for a suitable homogeneous ideal $J\subseteq P$.  The fact that $G$ is bipartite implies that the isomorphism above is graded (with respect to the standard grading on $P$).

\begin{prop}\label{RS}
Let $G$ be a bipartite graph satisfying the WSC property. Then there exists a term order $\prec$ on $P$ such that
\begin{compactenum}
\item $LT_{\prec}(J)$ is squarefree;
\item the simplicial complex $\Delta=\Delta(LT_{\prec}(J))$ is shellable.
\end{compactenum}
\end{prop}
\begin{proof}
The first item follows from the fact that $\AG$ is a normal semigroup ring (see Bruns-Gubeladze \cite[Corollary 7.20]{BG}).

As for the second item, by a result of Sturmfels (see \cite[Theorem 7.18 and Remark 7.21]{BG}) the simplicial complex associated to the radical of any initial ideal of a toric normal ideal is always shellable.
\end{proof}

\begin{cor}\label{coro}
Let $G$ be a bipartite graph satisfying the WSC property. Then for every $k>0$ \ $HF_{\AG}(k)=HP_{\AG}(k)$ (where $HP_{\AG}$ is the Hilbert polynomial of $\AG$);
\end{cor}
\begin{proof}
Choosing a term order $\prec$ as in Proposition \ref{RS} we have that, since $LT_{\prec}(J)$ is squarefree, $HF_{P/LT_{\prec}(J)}(k)=HP_{P/LT_{\prec}(J)}(k)$  for every $k>0$ (see \cite[Exercise 5.1.27]{BH}), so we conclude by noticing that $HF_{\AG}=HF_{P/LT_{\prec}(J)}$.
\end{proof}

Assume $G$ is unmixed without isolated points (by Remark \ref{thm:noisolated} there is no loss of generality in doing so). Consider the partition of the vertex set $[n]=A \cup B$. By the unmixed assumption, the basic 1-cover that yields 1 on $A$ and 0 on $B$ must have the same number of ones than the basic 1-cover that yields 1 on $B$ and 0 on $A$. Since $G$ has no isolated points, $|A|=|B|=\frac{n}{2}$.\\
We define a partial order on the set of basic 1-covers as follows:
\[\aa \leq \bb \iff a_i\leq b_i \ \ \forall \ i \in A.\]
It is quite easy to see that, by the unmixed assumption, the pair $\mathcal{L}=(\{\mbox{basic 1-covers}\},\leq)$ is a distributive lattice: in fact, given two basic 1-covers $\aa$ and $\bb$ one has only to check that the maximum $\aa \vee \bb$ and the minimum $\aa \wedge \bb$ are defined as follows:
\begin{displaymath}
(a \vee b)_i = \left\{ \begin{array}{cc} \max\{a_i,b_i\} & \mbox{if } i \in A\\
\min\{a_i,b_i\} &  \mbox{if } i\in B \end{array}  \right.
\end{displaymath}
\begin{displaymath}
(a \wedge b)_i = \left\{ \begin{array}{cc} \min\{a_i,b_i\} & \mbox{if } i \in A\\
\max\{a_i,b_i\} &  \mbox{if } i\in B \end{array}  \right.
\end{displaymath}
The fact that we have defined 1-covers is trivial, and that they are basic follows from the unmixed assumption. The distributivity of the two operations is straightforward.

Considering two incomparable basic 1-covers $\aa$ and $\bb$, it is clear that the binomial $X_{\aa} \cdot X_{\bb} - X_{\aa \vee \bb} \cdot X_{\aa \wedge \bb}$  belongs to $J$. Actually it is not difficult to see that the following equality holds:
\[(X_{\aa} \cdot X_{\bb} - X_{\aa \vee \bb}\cdot X_{\aa \wedge \bb})=J.\]
Consider a binomial $U-V$ belonging to $J$, but not in $(X_{\aa} \cdot X_{\bb} - X_{\aa \vee \bb}\cdot X_{\aa \wedge \bb})$. Suppose it is of minimal degree among those with this property. Using the relations $X_{\aa} \cdot X_{\bb} = X_{\aa \vee \bb} \cdot X_{\aa \wedge \bb}$ modulo $J$ one can assume that there exists a basic 1-cover $\cc$ (respectively $\dd$) maximal among the basic 1-covers ${\bf e}$ for which $X_{{\bf e}}$ divides $U$ (respectively $V$). It easily follows that $\cc=\dd$, and $J$ being prime this means that
\[\displaystyle \frac{U-V}{X_{\cc}}\in J \setminus (X_{\aa} \cdot X_{\bb} - X_{\aa \vee \bb}\cdot X_{\aa \wedge \bb}),\]
a contradiction to the minimality of the degree of $U-V$.

In particular we have proved the following (for the definition of Hibi ring see the paper of Hibi \cite{Hi})

\begin{thm}\label{hibiring}
Assume $G$ is bipartite unmixed with at least one edge. Then $\AG$ is a Hibi ring on $\mathcal{L}$ over $\KK$.
\end{thm}

Now choose a degrevlex term order $\prec$ on $P$ which is a linear extension of $\leq$ (i.e. $X_{\aa} \prec X_{\bb}$ provided $\aa < \bb$). Then
\[LT_{\prec}(J)=(X_{\aa} \cdot X_{\bb} : \aa \mbox{ and }\bb \mbox{ are not comparable}).\]
So with respect to this term order $LT_{\prec}(J)$ is the edge ideal $I(\mathcal{G})$ of the graph $\mathcal{G}$ on the vertex set $\{\mbox{basic 1-covers of }G\}$ whose edges are $\{\aa,\bb\}$ such that $\aa$ and $\bb$ are incomparable basic 1-covers. Moreover $\mathcal{G}$ is a Cohen-Macaulay graph again by \cite[Theorem 7.18,and Remark 7.21]{BG}.
Now recall some basic definitions for a poset $\Pi$:
\begin{compactitem}[--]
\item a chain of length $k$ in $\Pi$ is a sequence $\pi_k>\pi_{k-1}>\ldots >\pi_1$, \ $\pi_i \in \Pi$;
\item the rank of $\Pi$, denoted by $\rank(\Pi)$, is the maximum of the length of the chains in $\Pi$;
\end{compactitem}

\begin{prop}\label{unmixed}
Let $G$ be a bipartite unmixed graph with at least one edge. Then
\begin{compactenum}
\item $\AG$ is Koszul;
\item $\dim(\AG)= \rank(\mathcal{L})$;
\item the multiplicity $e(\AG)$ of $\AG$ is equal to the number of maximal chains in $\mathcal{L}$.
\end{compactenum}
\end{prop}
\begin{proof}
We can assume that $G$ has no isolated point by Remark \ref{thm:noisolated}.

(1). It follows from the fact, described above, that $J$ has a quadratic Gr\"obner basis.

(2). This is because $\AG$ is an ASL on $\mathcal{L}$ over $\KK$ by Theorem \ref{hibiring}.

(3). For the term order described above we have $LT_{\prec}(J)=I_{\Delta(\mathcal{L})}$ (see the book of Bruns and Herzog \cite[p. 208]{BH} for the definition of a simplicial complex associated to a poset). So the conclusion follows by the relation between the $h$-vector of $P/I_{\Delta(\mathcal{L})}$ and the $f$-vector of $\Delta(\mathcal{L})$ (\cite[Corollary 5.1.9]{BH}).
\end{proof}

\begin{os}
Point (2) of the above proposition is due to Bertone-Micale (\cite[Theorem 3]{BM}), which proved this result following another way.
\end{os}

We conclude this subsection giving a characterization of the Cohen-Macaulay bipartite graph for which $\AG$ is Gorenstein. Recall that a bipartite graph $G$ is Cohen-Macaulay if and only if $G=G(P)$ for a suitable poset $P$ (see \cite[Theorem 3.4]{HH}). Recall also that a poset $P$ is pure if all its maximal chains have the same length.

\begin{cor}\label{gor}
Assume $G=G(P)$ is Cohen-Macaulay. Then $\AG$ is Gorenstein if and only if $P$ is pure.
\end{cor}
\begin{proof}
In this case, one can easily check that $\mathcal{L}$ is the distributive lattice $\mathcal{J}(P)$ consisting of the poset ideals of $P$ ordered by inclusion. So the conclusion follows by \cite[main Corollary of p.105]{Hi}
\end{proof}

\subsection{Extension to Hypergraphs}\label{hyper}

Let us discuss a more general context, namely, hypergraphs. Our goal is to use algebraic tools to give upper and lower bounds for the dimension (Theorem\ref{generaldimension}).

A \textit{hypergraph} $H$ on $n$ vertices is a non-empty collection of subsets of $[n]$ (called ``faces''), such that if we take any two faces, the first is not contained in the second. A graph can be seen as a hypergraph where all faces have cardinality two. More generally, a simplicial complex together with the list of its facets forms a hypergraph; since all hypergraphs arise this way, in some sense the hypergraph notion is just a revisitation of the simplicial complex notion.

Let $H$ be a hypergraph and $F$ the set of the faces of $H$; let $\omega$ be a function that assigns to each face of $H$ a positive weight.
A non-zero vector $\aa \in \NN^n$ is called a $k$-cover (with respect to $\omega$) if
\[ \sum_{i \in f} a_i\geq k \omega(f) \ \ \ \forall \ f \in F. \]

We denote with $\AA_k \subseteq S[t]$ the $\KK$-vector space generated by all the
\[x_1^{a_1}\cdots x_n^{a_n}t^k\]
with $\aa$ a $k$-cover. The graded $S$-algebra
\[\AA=\bigoplus_{k\geq 0}\AA_k\]
is called the \textit{vertex cover algebra} of the hypergraph $(H,\omega)$ (see \cite{HHT}.)
One can then define the algebra
\[\bar{A}(H,\omega)=\AA /\mm \AA.\]
If we consider the $I$-filtration $\mathcal{F}=(I_{\omega,k})_{k \in \NN}$, where
\[I=\bigcap_{f\in F}\wp_f, \ \ \ I_{\omega,k}=\bigcap_{f \in F}\wp_f^{\omega(f)k} \hbox{ and } \wp_f=(x_i:i\in f),\]
it turns out that $\AA$ is the Rees algebra associated to $\mathcal{F}$, while $\bar{A}(H,\omega)$ is the fiber cone associated to the same filtration.
\[\bar{A}(H,\omega)=F_{\mm}(\mathcal{F})=\bigoplus_{k \geq 0}I_{\omega,k}/\mm I_{\omega,k}\]
(for the language of filtrations see \cite[Chapter 4, Section 5]{BH} or the manuscript of Rossi--Valla \cite{RV}).\\
In case $\omega$ is the constant function that assigns 1 to every vertex of $H$,  $A(H,\omega)$ (resp. $\bar{A}(H,\omega)$) is the symbolic Rees algebra (resp. the symbolic fiber cone) associated to the ideal $I$.

\begin{thm}\label{generaldimension}
Assume $\KK$ is an infinite field and let $M$ denote the maximal
cardinality of a face of $H$. Then
\[n - \depth(S/I) \leq \dim(\bar{A}(H,\omega)) \leq n - \left[ \frac{n-1}{M} \right].\]
\end{thm}
\begin{proof}
From \cite[Corollary 2.2]{HHT} there exists an integer $d$ such that
\[I_{\omega,d}^k=I_{\omega, dk} \ \ \ \forall \ k\geq 1.\]
This means that
\[ \bar{A}(H,\omega)^{(d)}=F_{\mm}(I_{\omega,d})=\bigoplus_{k \geq 0}I_{\omega,d}^k/\mm I_{\omega,d}^k, \]
where $\bar{A}(H,\omega)^{(d)}$ denotes the $d$-th Veronese subalgebra of $\bar{A}(H,\omega)$. Since the dimension of $\bar{A}(H,\omega)$ and the dimension of $\bar{A}(H,\omega)^{(d)}$ coincide, it suffices to prove the claim for $F_{\mm}(I_{\omega,d})$.\\
Set $\delta = \dim F_{\mm}(I_{\omega,d})$.

{\it The lower bound.} From \cite[Theorem 1 and 2]{NR} it follows that $\delta$ is the cardinality of a minimal set of homogeneous generators of a minimal reduction $J$ of $I_{\omega,d}$. So, since $\sqrt{J}=\sqrt{I_{\omega,d}}=I$, we have that $\operatorname{ara}(I) \leq \delta$. It is known that $\operatorname{cd}(S,I)\leq \operatorname{ara}(I)$, where $\operatorname{cd}(S,I)$ denotes the cohomological dimension of $I$, that is, the maximum integer $i$ for which the local cohomology module $H_I^i(S)\neq 0$.\\
Since $I$ is a squarefree monomial ideal, a theorem of Lyubeznik \cite{L1} implies that
\[\operatorname{cd}(S,I)=n - \depth(S/I),\]
whence the lower bound follows.

{\it The upper bound.} By a theorem of Burch (\cite{B}) \[ \delta \leq n-\min_{k\geq 1} \{ \depth(S/I_{\omega,d}^k) \}. \]
By our choice of $d$ we have that $I_{\omega,d}^k=I_{\omega,dk}$, so the associated prime ideals of $I_{\omega,d}^k$ are actually the minimal prime ideals of $I$ for each $k\geq 1$, i. e.
\[ \operatorname{Ass}(S/I_{\omega,d}^k)= \{\wp_f : f \in F\}. \]
Now we can use a result of Lyubeznik \cite[Proposition 2]{L2} and conclude that
\[ \min_{k\geq 1} \{ \depth(S/I_{\omega,d}^k)\} \geq \left[ \frac{n-1}{M} \right]. \]
\end{proof}

Notice that, if $G$ is a bipartite graph, by Theorem \ref{thm:inequalities} $\dim(\AG)\leq a+1\leq [n/2]+1$. The theorem above generalizes this inequality to every graph, not necessarily bipartite.

Next we give a combinatorial interpretation of Theorem \ref{generaldimension}.

\begin{cor}
Let $(H, \omega)$ be a weighted hypergraph on $n$ vertices, and $M$ the maximal cardinality of a face of $H$. The basic $k$-covers of $H$ with respect to $\omega$ are asymptotically counted (for $k$ large) by a polynomial $P$ in $k$, and
\[ M-1 \leq \deg P \leq (n - 1) - \left[ \frac{n-1}{M}\right]. \]
\end{cor}
\begin{proof}
The mentioned polynomial is the Hilbert polynomial of $\bar{A}(H, \omega)$, so its degree is equal to $\dim( \bar{A}(H,\omega))-1$.
The conclusion follows then from Theorem \ref{generaldimension}; for the left inequality we also used that $\depth(S/I)\leq n- \operatorname{ht}(\wp)$ for every $\wp \in \operatorname{Ass}(S/I)$.
\end{proof}

\subsection{Applications to the Arithmetical Rank}

In this subsection we relate the previous results to the arithmetical rank of certain squarefree monomial ideals. Recall the following result, which we already explained in the Definition section.


\begin{prop}\label{ara}
Let $G$ be a bipartite graph on $[n]$. Let $I$ be the Alexander dual of the edge ideal of $G$, that is,
\[ I = \bigcap_{\{i,j\} \hbox{edge of} G}(x_i,x_j). \]
If $\KK$ is an infinite field, \[\ara(I)\leq \dim(\bar{A}(G)).\]
\end{prop}

\begin{cor}\label{uara}
Let $G$ be a bipartite graph on $[n]$.
\begin{compactenum}
\item If $G$ is a tree, then $\ara(I)\leq \gdim(G)$.
\item If, in some standard drawing of $G$, the degree of the vertices below is at least $s$, then  $\ara(I) \leq a-s+2$.
\item If $G$ is unmixed, then $\ara(I)\leq \rank(\mathcal{L})$, where $\mathcal{L}$ is the distributive lattice associated to $G$.
\item If $G$ is a tree, then the Castelnuovo-Mumford regularity of the edge ideal of $G$ is at most $ \gdim(G)$.
\item If in some standard drawing for $G$ the degree of the vertices below is at least $s$, then the Castelnuovo-Mumford regularity of the edge ideal of $G$ is at most $ a-s+2$.
\item If $G$ is unmixed, then the Castelnuovo-Mumford regularity of the edge ideal of $G$ is at most $ \rank(\mathcal{L})$.
\end{compactenum}
\end{cor}

\begin{proof}
If $\KK$ is finite we can replace it by its algebraic closure $\bar{\KK}$: In fact, since $\KK \subseteq \bar{\KK}$ is a faithfully flat extension, $\ara(I)=\ara(I (S\otimes_{\KK}\bar{\KK}))$. Thus we can assume that $K$ is infinite.
The first two items follow then by Theorem \ref{6}, Theorem \ref{thm:7}, Proposition \ref{unmixed} and Proposition \ref{ara}.

The fourth, fifth and sixth item follow instead from the inequality
\[\projdim(S/I)=\cd(S,I)\leq \ara(I)\]
which the reader may find in \cite{L1}, and from Terai's identity \cite{Te} \[\operatorname{reg}(I^{\vee})=\projdim(S/I),\] where $I^{\vee}$ denotes the Alexander dual of $I$.
\end{proof}

Here are some examples:

\begin{es}\label{octagon}
Let $G$ be the octagon. Then
\begin{center}$I=(x_1x_3x_5x_7, \ \  x_2x_4x_6x_8, \ \ x_1x_3x_5x_6x_8, \ \ x_1x_3x_4x_6x_8, \ \ x_1x_3x_4x_6x_7,$ \\
$\qquad x_1x_2x_4x_5x_7, \ \ x_1x_2x_4x_6x_7, \ \ x_2x_3x_5x_6x_8, \ \ x_2x_3x_5x_7x_8, \ \ x_2x_4x_5x_7x_8).$
\end{center}
By the previous theorem $\ara(I)\leq 4$. Moreover by a Theorem of Lyubeznik \cite{L1}
\[\projdim(S/I)=\cd(S,I) \leq \ara(I),\]
and $\projdim(I)=4$ (using CoCoA \cite{Cocoa}), so
\[\ara(I)=4.\]
\end{es}

\begin{os}
In general, if $G=C_{2a}$ is the $2a$-cycle, then by Corollary \ref{uara} we obtain
\[ \ara(I)\leq a.\]
We know that equality holds if $a=2,3,4$; for the decagon, $\projdim(S/I)=4$, so we only know that
\[4 \leq \ara(I) \leq 5.\]
\end{os}

\begin{es}
For any natural numbers $2\leq r \leq a\leq b$ consider the graph of Example \ref{alex}. Then $\ara(I) \leq r+1$.\\
For example, if $r=4$, $a=b=6$ it is not difficult to see that the ideal is
\begin{center}
$I=(x_1x_2x_3x_5x_7x_9, \ \ x_1x_2x_3x_5x_7x_8x_{10}x_{11}x_{12}, \ \ x_1x_2x_3x_5x_6x_8x_{10}x_{11}x_{12}, \ \ x_1x_2x_3x_5x_6x_8x_9,$\\
$x_4x_6x_8x_{10}x_{11}x_{12}, \ \ x_4x_6x_8x_9, \ \ x_4x_6x_7x_9, \ \ x_4x_5x_7x_9, \ \ x_4x_5x_7x_8x_{10}x_{11}x_{12})\subseteq S=\KK[x_1, \ldots ,x_{12}].$
\end{center}
In this case $\ara(I)\leq 5$; however $\projdim(S/I)=4$, so $4\leq \ara(I)\leq 5$.
\end{es}

\begin{es}
Fix $a \geq 2$ and let $G=G_a$ be the $(a-1)$-regular bipartite graph on $2a$ vertices, already considered in Proposition \ref{min multiplicity}. Up to a change of coordinate, the Alexander dual of the edge ideal of $G_a$ is

\[I_a := \left(x_1 \cdots x_a ,  \ \ x_{a+1}\cdots x_{2a}, \ \ x_1 \cdots \widehat{x_i}\cdots x_a x_{a+1}\cdots \widehat{x_{a+i}}\cdots x_{2a} \; ~|~  \; i=1, \ldots ,a \right),\]
where $\widehat{x_j}$ means that we are skipping the variable $x_j$. By Theorem \ref{uara} $\ara(I_a)\leq 3$ for every $a$. We claim that equality holds.

To prove this we argue by induction on $a$, and prove that $\cd(\KK[x_1,\ldots, x_{2a}],I_a)=3$.\\

The case $a=2$ can be checked with CoCoA \cite{Cocoa}. Suppose now $a>2$, and specialize at $x_a = x_{2a} = 1$. Since local cohomology is independent of the base ring chosen,
\[ H_{I_a}^3(\KK[x_1, \ldots  ,\widehat{x_a}, \ldots, x_{2a-1}]) \cong H_{I_{a-1} }^3(\KK[x_1, \ldots , \widehat{x_a}, \ldots, x_{2a-1}]) \neq 0, \]
whence we conclude. In particular,
\[\ara(I_a)=3 \ \ \forall \ a\geq 2.\]
\end{es}

\subsection{Acknowledgements.}
The authors wish to thank Juergen Herzog and Volkmar Welker (together with all the organizers of Pragmatic 2008, Catania, Italy) for many helpful lessons, suggestions, and discussions.
Many thanks also to Aldo Conca for proofreading the paper.

\end{document}